\documentclass{amsart}

	\usepackage{amsmath}
	\usepackage{amssymb}
	\usepackage{graphicx}
	\usepackage[dvipsnames]{xcolor}
	\usepackage{amsxtra}
	\usepackage{mathtools}
	\usepackage{enumerate}
	\usepackage{enumitem}
	\usepackage{nicefrac}
	\mathtoolsset{showonlyrefs}
	\usepackage{mathrsfs}
	\usepackage{hyperref}
	\hypersetup{
	  colorlinks   = true,
	  urlcolor     = blue,
	  linkcolor    = blue,
	  citecolor   = ForestGreen
	}

	\newtheorem{theorem}{Theorem}[section]
	\newtheorem{lemma}{Lemma}[section]
	
	\newtheorem{corollary}[lemma]{Corollary}
	\newtheorem{conjecture}[lemma]{Conjecture}

	\theoremstyle{definition}
	
	\newtheorem{remark}[lemma]{Remark}

	\DeclareMathOperator{\spn}{span}
	\newcommand {\LL}{L_{a,b}}
	\newcommand{\R}{\mathbb{R}}
	\newcommand{\HH}{\mathcal{H}}
	\newcommand{\DD}{\mathcal{D}}

\title[Mixed local and nonlocal Lane--Emden equation]
{Mixed local and nonlocal Laplacian without standard critical exponent for Lane-Emden equation} 

\author[B. Barrios]{Bego\~na Barrios}
	\address{Bego\~na Barrios \hfill\break\indent
		Departamento de An\'{a}lisis Matem\'{a}tico e IMAULL
		\hfill \break \indent Universidad de La Laguna
		\hfill \break \indent C/. Astrof\'{\i}sico Francisco S\'{a}nchez s/n,
		38200 -- La Laguna, SPAIN}
	\email{bbarrios@ull.edu.es}

\author[L. M. Del Pezzo]{Leandro M. Del Pezzo}
	\address{Leandro M. Del Pezzo \hfill\break\indent
		IESTA --Facultad de Ciencias Económicas y de Administración
		\hfill\break\indent
		Universidad de la República\hfill\break\indent
		Av. Gonzalo Ramírez 1926, 11200 Montevideo,\hfill\break\indent
		Departamento de Montevideo - Uruguay.}
	\email{leandro.delpezzo@fcea.edu.uy}

\author[A. Quaas]{Alexander Quaas}
	\address{A. Quaas
		\hfill\break\indent Departamento de Matem\'atica,
		Universidad T\'ecnica Federico Santa Mar\'ia
		\hfill\break\indent Casilla V-110, Avda. Espa\~na, 1680 --
		Valpara\'iso, CHILE.}
	\email{alexander.quaas@usm.cl}

\begin{document}

	\begin{abstract}
		In this paper, we investigate a mixed elliptic equation
		involving both local and nonlocal Laplacian operators,
		with a power-type nonlinearity. Specifically, we consider a
		Lane--Emden type equation of the form
		\[
			-\Delta u + (-\Delta)^s u
				= u^p\quad\text{ in }\R^n,
		\]
		where the operator combines the classical Laplacian and
		the fractional Laplacian, \(s\in(0,1)\). We first observe, by a
		variational argument that works for any finite sum of fractional
		Laplacians of different orders, that positive radial solutions exist
		in the whole range \(\frac{n+2s}{n-2s}<p<\frac{n+2}{n-2}\), while in
		star-shaped bounded domains the threshold for existence is the local
		Sobolev exponent \(\frac{n+2}{n-2}\). Hence, for a sum of two diffusion
		operators without zero-order term, the classical duality between the
		critical exponent in \(\R^n\) and the one in star-shaped domains fails.
		Our main results describe the solutions near both thresholds. For
		\(n\ge6\) (or \(n=5\) and \(s>\frac14\)) and \(p<\frac{n+2}{n-2}\) close to
		\(\frac{n+2}{n-2}\), we construct by a Lyapunov--Schmidt reduction radial
		solutions which concentrate at the origin with the profile of the
		Aubin--Talenti bubble, at the explicit rate
		\(\delta_\varepsilon^{2(1-s)}\sim\varepsilon=\frac{n+2}{n-2}-p\) dictated by the
		competition between the two diffusions. Dually, in every dimension \(n\ge3\)
		and for \(p>\frac{n+2s}{n-2s}\) close to \(\frac{n+2s}{n-2s}\), we construct
		radial solutions which bifurcate from the bubble of \((-\Delta)^s\) and vanish
		by spreading out at the rate \(\mu_\eta^{-2(1-s)}\sim\eta=p-\frac{n+2s}{n-2s}\).
		These results show that this mixed Lane--Emden equation does not admit a
		critical exponent in the traditional sense.
	\end{abstract}
\maketitle
\section{Introduction}

	The objective of this paper is to explore critical exponents
	associated with a Lane--Emden type equation
	involving a mixed elliptic operator of the form
	\[
		\LL(u)\coloneqq -a\Delta u +b (-\Delta)^s u = u^p,
	\]
	with \(a,b\in [0,\infty),\) with at least one of them strictly
	positive, and \(s\in(0,1)\). As is well known, the fractional Laplacian is defined by
	the singular integral
	\[
		(-\Delta)^s u(x)\coloneqq C_{n,s}\lim_{r\to 0^+}
		\int_{\R^n\setminus B(x,r)}
		\frac{u(x)-u(y)}{|x-y|^{n+2s}}dy,
	\]
	where \(C_{n,s}>0\) is the usual normalization constant, chosen so that
	the Fourier symbol of \((-\Delta)^s\) is \(|\xi|^{2s}\), that is,
	\(\widehat{(-\Delta)^su}(\xi)=|\xi|^{2s}\widehat u(\xi)\). Accordingly, throughout the paper we
	use the bilinear form
	\[
		\langle u,v\rangle_s\coloneqq \frac{C_{n,s}}{2}\int_{\R^n}\int_{\R^n}
		\frac{(u(x)-u(y))(v(x)-v(y))}{|x-y|^{n+2s}}\,dy\,dx,\qquad
		[u]_s^2\coloneqq\langle u,u\rangle_s,
	\]
	so that \(\int_{\R^n}(-\Delta)^su\, v\,dx=\langle u,v\rangle_s\) for smooth,
	decaying \(u,v\).

	Mixed operators have attracted growing attention, mainly because they
	model situations where local diffusion (governed by Brownian motion) and
	nonlocal, L\'evy-type diffusion coexist (see, e.g., \cite{BonforteVa}).
	Such combinations show up in materials with heterogeneous microstructure,
	in financial models with both local volatility and jumps, and in
	population dynamics where individuals move through both short- and
	long-range interactions. They are also mathematically interesting in
	their own right, since the local and nonlocal parts compete in
	regularity and in scaling. A non-exhaustive list of references on the
	topic is \cite{Barles1, Barles2,Aba, Bia, Bia1,Bia2,Biswas,fil,Ga}, see
	also \cite{DSVZ,FangShangZhang,SVWZ}.
	
	\medskip

	Concerning critical exponents, recall that the critical exponent in the whole space associated with
	\[
		F(u)= u^p,
	\]
	where \(F\) is a (possibly nonlinear) elliptic operator,
	local or nonlocal in nature , without zero-order term, is defined by
	\[
		p_*(F)\coloneqq \inf \{p>1\colon
			\text{{\small there is a positive
				solution to} } F(u)=u^p \text{ {\small in }}
					\R^n \}.
	\]
	For operators with a zeroth-order term, such as \(-\Delta u+u\), constant
	solutions and the different behaviour at infinity make this notion ambiguous, see
	Remark~\ref{rem:zero-order}. For this reason we restrict ourselves to sums of
	homogeneous diffusion operators.
	
	The critical exponent for Lane--Emden equations runs through much of
	PDE theory: in the classical case \(F=-\Delta\) it connects the Yamabe
	problem with the sharp constants in Sobolev's inequality. See for
	instance \cite{Yamabe,Talenti,Lieb}.
	
	\medskip

	On the other hand, the critical exponent in a bounded
	(smooth) domain is defined by
	\[
		p^*(\Omega, F)\coloneqq
			\sup \left\{p>1\colon
			\begin{aligned}
			&\text{{\small there is a
			positive solution to} } \\
			&F(u)= u^p  \text{ {\small in} } \Omega, u=0\text{ {\small on} } \Gamma
			\end{aligned}
			\right\}.
	\]
	Here \(\Gamma=\partial\Omega\) if \(F\) is a local operator and
	\(\Gamma=\R^n\setminus\Omega\) if \(F\) is a  nonlocal operator.

	For the classical Laplacian \(F=-\Delta\)  and when \(\Omega\)
	is star-shaped,
	it is well known that a form of {\it duality}
	appears, that is,
	\[
		p^*(\Omega, -\Delta)=p_*(-\Delta)=p_1\coloneqq\frac{n+2}{n-2},
	\]
	see \cite{Poho, Caff1, Chen, Gidas}.
	In the radial case,
	this duality of the exponent (that is,
	\(p^*(B_R, -\Delta)=p_*(-\Delta)\))  is clear,
	since a solution \(u_\gamma(r)\), \(r \in (0,+\infty)\), is a
	solution of an ODE with positive initial condition
	\(u(0)=\gamma>0\)
	that crosses or not the  \(r\)-axis. In other words,  if
	\(1<p<p_1\) then there exists \(r_0\) such that \(u_\gamma(r_0)=0,\)
	and if \(p\ge p_1\) then \(u_\gamma(r)>0\) for all
	\(r \in (0,+\infty).\) By scaling, \(\gamma\) can
	be an arbitrary value of \(\R^+\). The Emden--Fowler transformation is a classical change of variable
	to study these types of equations and gives the radial duality
	described above.

	\medskip

	It is then natural to look for domains
	\(\Omega\) with non-trivial topology (not star-shaped) such that
	\(p^*(\Omega, -\Delta)>p_1\). See for instance \cite{manuel} and also the starting
	point in \cite{Bahri} at the critical exponent \(p_1\).  The case of the exterior domain is discussed for example in
	\cite{Davila} and the case of half space in \cite{Souplet}.
	Such domains are also of interest for mixed operators, but they fall
	outside the scope of this work, which focuses on the whole space.

	\medskip

	Returning to star-shaped domains or the whole space,
	similar results to the local ones hold when considering the purely
	nonlocal operator
	\(F=L_{0,1}=(-\Delta)^s\) (i.e. \(a=0\) and \(b=1\)).
	In fact, if \(\Omega\) is star-shaped, then we also
	have duality of the exponents and
	\[
		p^*(\Omega, (-\Delta)^s)=p_*((-\Delta)^s)
		=p_s\coloneqq\frac{n+2s}{n-2s},
	\]
	see for instance \cite{ros,val,ChenLiOu}. 
	
	Variational methods guarantee existence of solutions
	in any regular bounded domain \(\Omega\) for
	the subcritical range \(1<p<p_s,\) when \(\Omega\) is
	not necessarily star-shaped, see \cite{val}. Results on the critical case in the purely nonlocal setting can be found,
	for example, in
	\cite{bego}, where a H\'enon-type equation is discussed.
	
	\medskip

	Other equations exhibiting a coexistence of sub- and supercritical
	behaviour similar to that of \eqref{elproblema} are, for instance,
	\[
		\Delta u+u^p+u^q=0\quad\text{in }\quad\R^n,\, 1<p<p_1<q,
	\]
	studied in the seminal work \cite{Ni}. For that equation, existence of positive
	solutions with fast decay is shown in \cite{Bamon}, see also \cite{Campos}.
	Another closely connected problem with
	sub and supercritical
	behavior is the Matukuma type problem of the form
	\[
		\Delta u +k(r)u^p=0\, \text{ in }\R^n,
	\]
	studied in \cite{Felmer}, where fast decay solutions
	exist for special \(k\)
	having sub and supercritical behavior.
	Moreover,  the case where \(k\) is a mixed H\'enon type
	with sub and supercritical behavior is studied in
	\cite{Salo}. Fast decay solutions also exist.

	\medskip

	Returning to the case of mixed operators,
	if \(\Omega\) is a star-shaped domain and \(a,b>0\),
	it has been shown in Theorem 1.3 in \cite{ros2} (see also Example iv
	in Section 2) that \(p^*(\Omega, \LL) =p_1\). This value of the critical exponent also
	appears in a Brezis--Nirenberg type result \cite{BDVV-BN}. Using
	variational methods, existence results in bounded domains
	have been established for the subcritical range
	\(1<p<p_1\) in \cite{varational1}, even when \(b<0\).

	\medskip

	This leads to the  natural question of whether, for star-shaped domains,
	\[
		p_*(\LL)=p^*(\Omega, \LL)=p_1,
	\]
	that is, whether the duality also holds for
	the mixed operator.

	\subsection*{Existence in the whole intermediate range}
	The first observation of this paper is that the answer is negative,
	and that this follows from a direct variational argument. \footnote{we thank one
	of the referees for pointing out the argument in the case of two operators.}
	In fact the argument works for an arbitrary finite sum of
	fractional Laplacians of different orders. Let \(m\ge2\),
	\begin{equation}\label{eq:orders}
		0<s_1<s_2<\dots<s_m\le 1,\qquad c_1,\dots,c_m>0,
	\end{equation}
	where, for \(s_m=1\), \((-\Delta)^{1}=-\Delta\), and consider
	\begin{equation}\label{eq:suma}
		\mathcal{L}u\coloneqq\sum_{i=1}^m c_i(-\Delta)^{s_i}u=u^p
		\quad\text{in }\R^n,\qquad n\ge 3.
	\end{equation}
	We write \(p_{\kappa}\coloneqq\frac{n+2\kappa}{n-2\kappa}\) and
	\(2^*_\kappa\coloneqq p_\kappa+1=\frac{2n}{n-2\kappa}\) for \(\kappa\in(0,1],\)
	in particular \(p_1=\frac{n+2}{n-2}\) and \(2^*\coloneqq 2^*_1\).

	\begin{theorem}\label{thm:variational}
		Assume \eqref{eq:orders} and \(p_{s_1}<p<p_{s_m}\). Then \eqref{eq:suma}
		has a positive, radially symmetric and radially nonincreasing weak
		solution \(u\) with finite energy, that is,
		\(\sum_{i=1}^m c_i[u]_{s_i}^2<\infty\) (where \([u]_1=\|\nabla u\|_{L^2}\)).
		The solution is obtained, up to a multiplicative constant, as a minimizer of
		\begin{equation}\label{eq:Sq}
			S_q\coloneqq\inf\left\{\sum_{i=1}^m c_i[u]_{s_i}^2\colon
			u\in\DD,\ \|u\|_{L^q}=1\right\},\qquad q=p+1,
		\end{equation}
		where \(\DD\) is the energy space defined in Section~\ref{sec:variational}.
	\end{theorem}

	In particular, taking \(m=2\), \(s_1=s\), \(s_2=1\), for every \(a,b>0\) the
	equation \(\LL u=u^p\) has a positive radial solution in \(\R^n\) for every
	\(p_s<p<p_1\). Together with the nonexistence result \cite[Theorem 1.3]{ros2}
	for \(p\ge p_1\) and the existence of solutions in bounded domains for
	\(1<p<p_1\) \cite{varational1}, this yields:

	\begin{corollary}\label{cor:duality}
		Let \(a,b>0\), \(s\in(0,1)\), \(n\ge3\) and let \(\Omega\) be a bounded smooth
		star-shaped domain. Then
		\[
			p_*(\LL)\le p_s<p_1=p^*(\Omega,\LL).
		\]
	\end{corollary}

	Thus, for sums of homogeneous diffusion operators of different orders (with
	no zero-order term), the duality between the two critical exponents fails.
	We are not aware of previous results showing this phenomenon, note that
	equations with mixed operators and power nonlinearities have mostly been
	studied in the presence of zero-order terms or under mass constraints, see,
	e.g., \cite{DSVZ} and the references therein.

	\begin{remark}\label{rem:zero-order}
		The restriction to operators without zero-order terms is essential
		for the comparison above to be meaningful. For instance, positive
		finite-energy solutions of \(-\Delta u+u=u^p\) in \(\R^n\) exist in the whole
		range \(1<p<p_1\), whereas the Dirichlet threshold in star-shaped domains is
		\(p_1\). Moreover, constant functions solve this equation for every \(p\).
	\end{remark}

	\begin{remark}[Pohozaev identity and finite-energy solutions]\label{rem:Poho}
		The range of Theorem~\ref{thm:variational} is sharp for finite-energy
		solutions satisfying the Pohozaev identity. Indeed, if \(u\) is a positive solution of \eqref{eq:suma} with
		enough regularity and decay to justify the Pohozaev identity in \(\R^n\)
		(obtained by the usual dilation argument, see \cite{ros} for the nonlocal
		terms), then, writing \(M_i\coloneqq c_i[u]_{s_i}^2>0\),
		\[
			\sum_{i=1}^m M_i=\int_{\R^n}u^{p+1}\,dx,\qquad
			\sum_{i=1}^m\frac{n-2s_i}{2}M_i=\frac{n}{p+1}\int_{\R^n}u^{p+1}\,dx.
		\]
		If \(p+1\le 2^*_{s_1}\), then \(\frac{n}{p+1}\ge\frac{n-2s_1}{2}\ge\frac{n-2s_i}{2}\)
		with strict inequality for \(i\ge2\), which contradicts the two identities.
		Whereas if \(p+1\ge 2^*_{s_m}\) one argues in the same way with the reverse
		inequalities. Hence, among finite-energy solutions for which the Pohozaev identity holds,
		solutions can exist only for \(p_{s_1}<p<p_{s_m}\), and Theorem~\ref{thm:variational}
		shows that they do exist in this whole range. In particular, unlike in the purely local case, \(\LL\) with \(a,b>0\) has no such
		solutions at \(p=p_1\).
	\end{remark}

	Remark~\ref{rem:Poho} only concerns solutions with finite energy (or with
	enough integrability). This leads to the following conjecture, which remains open.

	\medskip
	\begin{conjecture}\label{conj}
		\(p_*(\LL)=p_s.\) More generally, \(p_*(\mathcal{L})=p_{s_1}\) for the operator
		in \eqref{eq:suma}. Equivalently, in view of Theorem~\ref{thm:variational},
		there are no positive solutions of \(\mathcal{L}u=u^p\) in \(\R^n\) for
		\(1<p<p_{s_1}\), without any integrability assumption on \(u\).
	\end{conjecture}
	\medskip

	A recent result supports this conjecture in the case \(m=2\), \(s_1=s\), \(s_2=1\): for
	the single mixed operator \(-\Delta u+(-\Delta)^su\), \cite{GuoLiXie} proves that positive
	distributional supersolutions of \(-\Delta u+(-\Delta)^su=u^q\) in \(\R^n\) exist if and only if
	\(q>\frac n{n-2s}\), the Serrin-type exponent dictated by the lower-order operator
	\((-\Delta)^s\). Since \(\frac n{n-2s}<p_s=\frac{n+2s}{n-2s}\), this does not reach
	Conjecture~\ref{conj}, which asks for nonexistence all the way up to \(p_s\). Still, it does
	confirm, for the smaller range \(1<p\le\frac n{n-2s}\), that no positive solution exists even
	without any energy or sign assumption.

	\subsection*{Concentrating solutions}
	Apart from radial symmetry and monotonicity, the variational argument of
	Theorem~\ref{thm:variational} gives no information on the shape of the
	solutions. The main results of this paper
	give a precise qualitative description of solutions near both thresholds.
	Near \(p_1\) we construct solutions which bifurcate from the Aubin--Talenti
	bubble of the local operator and concentrate at the origin, and near \(p_s\), solutions
	which bifurcate from the bubble of \((-\Delta)^s\) and spread out. In both
	cases the profile is explicit, and the rate is dictated by the competition
	between the two diffusions.

	Existence of positive solutions typically changes abruptly at the
	critical exponent, and locating this threshold for operators such as
	\(-\Delta+(-\Delta)^s\), where local and nonlocal effects compete, is
	exactly the question we address here.

	More precisely, we want to establish the existence
	of radial solutions to
	\begin{equation}
		\label{elproblema}
		\begin{cases}
			-\Delta u+ (-\Delta)^s u= n(n-2)u^{p_\varepsilon}
			&\text{in }\R^n,\, p_{\varepsilon}\coloneqq p_1-\varepsilon,\, n\ge5,\ n>6-4s,\\
			u>0  &\text{in }\R^n,
		\end{cases}
	\end{equation}
	where, for simplicity and
	without loss of generality, we have taken
	\[
		a=b=\frac{1}{n(n-2)}
	\]
	in \(L_{a,b}\), and  \(\varepsilon>0\)
	will be a small parameter. It is clear that  \(u_\varepsilon\) is a solution of
	\eqref{elproblema} if and only if
	\[
		u_{\delta, \varepsilon}(x)\coloneqq \delta^{\frac2{p_\varepsilon-1}}
		u_\varepsilon(\delta x),
		\text{ with }\delta>0,
	\]
	is a solution of
	\begin{equation}
		\label{elproblemaDelta}
		\begin{cases}
			-\Delta u+ \delta^{2(1-s)}(-\Delta)^s u=
			n(n-2)u^{p_\varepsilon} &\text{in }\R^n,\\
			u>0  &\text{in }\R^n,
		\end{cases}
	\end{equation}
	in the energy space \(\HH\) defined in \eqref{elespacio} below. As
	\(\delta\to0\), \eqref{elproblemaDelta} formally becomes
	the local problem
	\begin{equation}
		\label{elproblemalocalconp0}
		\begin{cases}
			-\Delta u= n(n-2)u^{p_1} &\text{in }\R^n,\\
			u>0  &\text{in }\R^n,
		\end{cases}
	\end{equation}
	whose finite-energy radial solutions are exactly the Talenti functions
	\begin{align}\label{bubble}
		U_{\delta} (x)\coloneqq\delta^{-\frac{n-2}2}
		U\left(\frac{x}{\delta}\right),\quad U(x)\coloneqq \frac1{(1+|x|^2)^{\frac{n-2}2}},\quad \delta>0,
	\end{align}
	all positive finite-energy solutions are the translates \(U_\delta(x-y)\), \(y\in\R^n,\)  
	see \cite{Talenti}.
	
	\medskip
	
	Our main result is the following.

	\begin{theorem}	\label{ElTeorema}
	 	 Let \(s\in(0,1)\) and \(n\ge5\) with \(n>6-4s\) (that is, \(n\ge6\), or \(n=5\) and
	 	 \(s>\frac14\)). There exist \(\varepsilon_0>0\)
	 	 and, for every \(\varepsilon\in(0,\varepsilon_0)\), a number
	 	 \(\delta_\varepsilon>0\) with
	 	 \[
	 	 	\delta_\varepsilon^{2(1-s)}=\lambda_\varepsilon\,\varepsilon,\qquad
	 	 	\lambda_\varepsilon\to\lambda_0\coloneqq
	 	 	\frac{(n-2)^3}{4(1-s)}\,\frac{\int_{\R^n}U^{2^*}dx}{[U]_s^2}\quad
	 	 	\text{as }\varepsilon\to0,
	 	 \]
	 	 (in particular \(\delta_\varepsilon\to0\) as \(\varepsilon\to0\)), such that problem
	 	 \eqref{elproblemaDelta} with \(\delta=\delta_\varepsilon\) has a radial solution of the form
		\[
			z_\varepsilon(x)= U(x) +\phi_\varepsilon(x),
			\qquad \|\phi_\varepsilon\|_{\mathcal{D}^{1,2}}
			+\delta_\varepsilon^{1-s}[\phi_\varepsilon]_s\le C\varepsilon.
		\]
		Consequently, the function
		\[
			u_\varepsilon(y)\coloneqq\delta_\varepsilon^{-\frac{2}{p_\varepsilon-1}}
			z_\varepsilon\Big(\frac{y}{\delta_\varepsilon}\Big)
		\]
		is a positive radial solution of \eqref{elproblema} which concentrates at
		the origin as \(\varepsilon\to0\): \(\delta_\varepsilon^{\frac{2}{p_\varepsilon-1}}
		u_\varepsilon(\delta_\varepsilon\,\cdot)\to U\) in \(\mathcal{D}^{1,2}(\R^n)\) and
		\[
			|\nabla u_\varepsilon|^2dx\rightharpoonup
			\Big(\int_{\R^n}|\nabla U|^2dx\Big)\,\delta_0
			\quad\text{weakly as measures.}
		\]
	\end{theorem}

	\begin{remark}Two remarks on this threshold.
		\begin{enumerate}[label=\textup{(\alph*)}]
			\item The value of \(\lambda_0\) is \(-A/B\), where \(A\) and \(B\) are given in
			Lemma~\ref{taqui}, it reflects the balance between the subcritical
			perturbation \(\varepsilon\) of the power, which tends to make the bubble
			spread and the lower-order diffusion \((-\Delta)^s\), which tends to
			make it shrink.
			
			\item The condition \(n>6-4s\) is used only to have
			\((-\Delta)^sU\in L^{2n/(n+2)}(\R^n)\) (Lemma~\ref{lem:fracU}), and we expect it to
			be technical, together with \(n\ge5\) it only excludes \(n=5\),
			\(s\le\frac14\). The assumption \(n\ge5\) is used separately (to have
			\(U,\psi\in H^1(\R^n)\)), so \(n=4\) is not covered either. For \(n=3\) and
			\(s\le\frac12\) there is a genuine obstruction, since then \([U]_s=\infty\).
		\end{enumerate}
	\end{remark}

	\subsection*{Blow-down at the lower exponent}
	The roles of the two operators can be exchanged. Let
	\[
		U_s(x)\coloneqq\frac{1}{(1+|x|^2)^{\frac{n-2s}{2}}},
	\]
	which satisfies \((-\Delta)^sU_s=c_{n,s}U_s^{p_s}\) in \(\R^n\) with
	\(c_{n,s}=2^{2s}\Gamma\left(\frac{n+2s}{2}\right)/\Gamma\left(\frac{n-2s}{2}\right)\) (for \(s=1\) this
	is \(n(n-2)\)), and whose dilations \(\mu^{-\frac{n-2s}{2}}U_s(x/\mu)\) are, up to
	translations, the positive solutions of this equation \cite{ChenLiOu,Lieb}.
	For \(\eta>0\) small we consider
	\begin{equation}\label{elproblema-s}
		\begin{cases}
			-\Delta u+(-\Delta)^su=c_{n,s}\,u^{p_s+\eta}&\text{in }\R^n,\\
			u>0&\text{in }\R^n.
		\end{cases}
	\end{equation}
	If \(z\) solves
	\begin{equation}\label{elproblemaSigma}
		(-\Delta)^sz+\sigma(-\Delta)z=c_{n,s}\,z^{p_s+\eta}\quad\text{in }\R^n,\qquad z>0,
	\end{equation}
	with \(\sigma=\mu^{-(2-2s)}\), then \(u(x)\coloneqq\mu^{-\frac{2s}{p_s+\eta-1}}z(x/\mu)\) solves
	\eqref{elproblema-s}. As \(\sigma\to0\), \eqref{elproblemaSigma} formally becomes
	the critical fractional problem, and the relevant regime is now \(\mu\to\infty\).

	\begin{theorem}\label{TeoremaBlowDown}
		Let \(s\in(0,1)\) and \(n\ge3\). There exist \(\eta_0>0\) and, for every
		\(\eta\in(0,\eta_0)\), a number \(\sigma_\eta>0\) with
		\[
			\sigma_\eta=\lambda_\eta\,\eta,\qquad \lambda_\eta\to\lambda_s\coloneqq
			\dfrac{c_{n,s}(n-2s)^2}{4n(1-s)}\,\dfrac{\displaystyle\int_{\R^n}U_s^{2^*_s}dx}{\|\nabla U_s\|_{L^2}^2}
			\quad\text{as }\eta\to0,
		\]
		such that \eqref{elproblemaSigma} with \(\sigma=\sigma_\eta\) has a radial solution
		\(z_\eta=U_s+\phi_\eta\) with 
		\[
			[\phi_\eta]_s+\sigma_\eta^{1/2}\|\nabla\phi_\eta\|_{L^2}\le C\eta.
		\]
		Consequently, with \(\mu_\eta\coloneqq\sigma_\eta^{-\frac{1}{2-2s}}\to\infty\), the function
		\[
			u_\eta(x)\coloneqq\mu_\eta^{-\frac{2s}{p_s+\eta-1}}z_\eta\Big(\frac{x}{\mu_\eta}\Big)
		\]
		is a positive radial solution of \eqref{elproblema-s} which vanishes by
		spreading out as \(\eta\to0\): \(\int_{\R^n}u_\eta^{2^*_s}dx\to\int_{\R^n}U_s^{2^*_s}dx\), while
		\(\int_{B_R}u_\eta^{2^*_s}dx\to0\) for every \(R>0\).
	\end{theorem}

	Theorems~\ref{ElTeorema} and~\ref{TeoremaBlowDown} are an exact counterpart of
	each other, in two respects. Quantitatively, \(\delta_\varepsilon\sim(p_1-p)^{\frac{1}{2-2s}}\to0\)
	while \(\mu_\eta\sim(p-p_s)^{-\frac{1}{2-2s}}\to\infty\). The two statements are not symmetric in
	generality, though, since Theorem~\ref{TeoremaBlowDown} holds in every dimension \(n\ge3\)
	(the error term \(\Delta U_s\) belongs to \(L^{2n/(n+2s)}(\R^n)\) for every such \(n\)), whereas
	Theorem~\ref{ElTeorema} needs \(n>6-4s\), coming from \((-\Delta)^sU\in L^{2n/(n+2)}(\R^n)\).
	Qualitatively, at \(p_1\) the local bubble is selected and the lower-order term \((-\Delta)^s\)
	forces it to shrink, the reduced equation reading \(\varepsilon A+\tau B+o(\varepsilon)=0\) with
	\(A>0>B\). At \(p_s\) the fractional bubble is selected and the higher-order term \(-\Delta\)
	forces it to spread, with \(-\eta A_s+\sigma B_s+o(\eta)=0\) and \(A_s,B_s>0\). In the second case
	the solutions lose all their \(L^{2^*_s}\) mass on compact sets as \(p\downarrow p_s\), consistent
	with Conjecture~\ref{conj}.

	To prove Theorem~\ref{ElTeorema}, we use a
	Lyapunov--Schmidt reduction method, which was first used in \cite{Floer} (see also \cite{NiTakagi}) in the
	context of partial differential equations and has become a very
	active area in PDEs since then.
	For global bifurcation results for semilinear equations in the whole space
	see \cite{BrownStavrakakis}. For a review of some other perturbation methods
	in \(\R^n\) see the book \cite{Mal} and also
	\cite{Campos,Salo}. In our particular case we follow the
	Lyapunov--Schmidt reduction in the spirit of \cite{MR1938656}, where the
	main idea is to solve
	the problem in the orthogonal complement of the
	kernel of the linearized operator (a finite-dimensional reduction) and
	then to choose the parameter \(\delta\) so that the remaining component,
	along the kernel, also vanishes. The new difficulty is that the
	perturbation \(\delta^{2(1-s)}(-\Delta)^s\) is singular at infinity: it
	dominates the Laplacian at scales \(|x|\gtrsim\delta^{-1}\). We therefore work in
	the \(\delta\)-dependent energy space \(\HH\) with the norm
	\(\|\nabla u\|_{L^2}^2+\delta^{2(1-s)}[u]_s^2\) and prove linear estimates that are
	uniform in \(\delta\). The fractional part of the norm is precisely what
	controls the slightly subcritical nonlinearity at infinity
	(Lemma~\ref{lem:interp}), and this works because the loss in the constants
	is of order \(\delta^{-O(\varepsilon)}\), which stays bounded in the regime
	\(\delta^{2(1-s)}\sim\varepsilon\).

	\subsection*{Organization of the paper}
	In Section~\ref{sec:variational} we
	prove Theorem~\ref{thm:variational}. In Section~\ref{section.psn} we introduce
	the functional framework for the reduction, prove the uniform invertibility
	of the projected linearized operator, and compute the signs of the two
	quantities that drive the reduction (Lemma~\ref{taqui}), in particular
	\(\langle U,\psi\rangle_s<0\), where \(\psi\) spans the radial kernel.
	Section~\ref{sec:existence} contains the reduction and the proof of
	Theorem~\ref{ElTeorema}. In Section~\ref{sec:blowdown} we prove
	Theorem~\ref{TeoremaBlowDown}, following the same scheme with the roles of
	the two operators exchanged.


\section{Proof of Theorem~\ref{thm:variational}}\label{sec:variational}

	Throughout this section we assume \eqref{eq:orders}, \(n\ge3\), and
	\(p_{s_1}<p<p_{s_m}.\) We set 
	\[
		q\coloneqq p+1\in(2^*_{s_1},2^*_{s_m}).
	\] 
	We use the
	convention \([u]_1=\|\nabla u\|_{L^2}\),
	\(\langle u,v\rangle_1=\int_{\R^n}\nabla u\cdot\nabla v dx\).
	Let
	\[
		\DD\coloneqq\Big\{u\in L^{2^*_{s_m}}(\R^n)\colon [u]_{s_i}<\infty
		\text{ for } i=1,\dots,m\Big\},\]
	\[
		\|u\|_\DD^2\coloneqq E(u)\coloneqq\sum_{i=1}^mc_i[u]_{s_i}^2 .
	\]
	
	The space \(\DD\) is a Hilbert space. Since its elements vanish at infinity (they belong
	to \(L^{2^*_{s_m}}\)), the homogeneous Sobolev inequalities
	\(\|u\|_{L^{2^*_{s_i}}}\le S_{s_i}[u]_{s_i}\) hold, and in particular
	\(\DD\hookrightarrow L^{2^*_{s_1}}(\R^n)\cap L^{2^*_{s_m}}(\R^n)\). By
	interpolation,
	\begin{equation}\label{eq:interp-sum}
		\|u\|_{L^q}\le\|u\|_{L^{2^*_{s_1}}}^{\vartheta}\|u\|_{L^{2^*_{s_m}}}^{1-\vartheta}
		\le C\|u\|_\DD,\qquad \frac1q=\frac{\vartheta}{2^*_{s_1}}+\frac{1-\vartheta}{2^*_{s_m}},
	\end{equation}
	so the constraint in \eqref{eq:Sq} is a \(C^1\) functional on \(\DD\) and
	\(S_q\ge C^{-2}>0\).

	\medskip

	We split the proof into four steps: a symmetric minimizing sequence,
	compactness, the Euler--Lagrange equation, and positivity.

	\medskip

	\noindent\emph{Step 1: Symmetric minimizing sequence.}
	For \(u\in\DD\), it is clear that $E(|u|)\le E(u)$
	and that the Schwarz symmetrization \(u^*\) of \(|u|\)
	satisfies \(\|u^*\|_{L^r}=\|u\|_{L^r}\) for every \(r\). Moreover, since by the P\'olya--Szeg\H{o} inequality for \(s_i=1\) and by \cite{AlmgrenLieb} for \(s_i\in(0,1),\)
	\[
		[u^*]_{s_i}\le[u]_{s_i}  \quad \forall i,
	\] 	
	there is a minimizing sequence \((u_k)\subset\DD\) of nonnegative, radially symmetric,
	nonincreasing functions with \(\|u_k\|_{L^q}=1\) and \(E(u_k)\to S_q\) (recall
	\eqref{eq:Sq}) as \(k\to\infty.\)
	
	\medskip

	\noindent\emph{Step 2: Compactness.}
	The sequence is bounded in \(\DD\), so, up to a subsequence, 
	\[
		u_k\rightharpoonup u\text{ weakly in } \DD\quad \text{ and }\quad  u_k\to u \text{ a.e.}
	\]
	Moreover, \(u\) is nonnegative, radial and nonincreasing. For every \(R>0\), \((u_k)\) is bounded in \(H^{s_m}(B_R)\)
	(the \(L^2(B_R)\) norm is controlled by the \(L^{2^*_{s_m}}\) norm), and since
	\(q<2^*_{s_m}\) the Rellich--Kondrachov theorem gives 
	\[
		u_k\to u \text{ in } L^q(B_R).
	\]
	At infinity we use monotonicity: for \(r=|x|\),
	\[
		u_k(r)^{2^*_{s_1}}|B_r|\le\int_{B_r}u_k^{2^*_{s_1}} dx\le C
		\quad\Longrightarrow\quad u_k(r)\le C r^{-n/2^*_{s_1}},
	\]
	and therefore
	\[
		\int_{|x|>R}u_k^q\,dx\le\Big(\sup_{|x|>R}u_k\Big)^{q-2^*_{s_1}}
		\int_{\R^n}u_k^{2^*_{s_1}}dx\le C R^{-\frac{n(q-2^*_{s_1})}{2^*_{s_1}}},
	\]
	uniformly in \(k.\) The same bound holds for \(u\). 
	The pointwise bound above is
	an elementary radial lemma of the type introduced by Strauss in \cite{Strauss}. Note that
	\(\tfrac{n}{2^*_{s_1}}=\frac{n-2s_1}{2}\), so it has the same decay as the fractional Strauss
	inequality 
	\[
		|u(x)|\le C|x|^{-\frac{n-2s}{2}}[u]_s,
	\] 
	for radial \(u\in\dot H^s(\R^n)\), which holds without monotonicity when \(\frac12<s<\frac n2,\) 
	see \cite{ChoOzawa} and \cite[Lemma 3.1]{SireWei}. Here monotonicity allows us to treat
	every \(s_1\in(0,1)\). Hence \(u_k\to u\) in \(L^q(\R^n)\),
	\(\|u\|_{L^q}=1\) and, by weak lower semicontinuity, \(E(u)\le S_q\). Thus \(u\) is a
	minimizer.
	
	\medskip

	\noindent\emph{Step 3: Euler--Lagrange equation.}
	By the Lagrange multiplier rule,
	\[
		\sum_{i=1}^m c_i\langle u,v\rangle_{s_i}=S_q\int_{\R^n}u^{q-1}v\,dx
		\qquad\forall v\in\DD,
	\]
	and \(w\coloneqq S_q^{1/(p-1)}u\) is a nonnegative, radial, nonincreasing,
	nontrivial weak solution of \eqref{eq:suma}.
	
	\medskip

	\noindent\emph{Step 4: Positivity.}
	Since \(w\) is radially nonincreasing, if \(w\) were not positive a.e. there would
	exist \(r_0>0\) with \(w=0\) a.e. in \(\{|x|>r_0\}\). Then take a nontrivial nonnegative
	function \(v\in C^\infty_c(\{|x|>r_0\})\). Thus
	\[
		\langle w, v\rangle_1=0\quad\text{ and } \int_{\R^n} w^pv \, dx=0,
	\]
	whereas, for \(s_i<1\), since \(w\,v\equiv0\),
	\[
		\langle w,v\rangle_{s_i}=-C_{n,s_i}\int_{\R^n}\int_{\R^n}
		\frac{w(y)v(x)}{|x-y|^{n+2s_i}}\,dy\,dx<0 .
	\]
	Since \(m\ge2\), at least one \(s_i\in(0,1)\), and testing the equation
	with \(v\) gives a contradiction. Hence \(w>0\) a.e., which proves
	Theorem~\ref{thm:variational}.\qed

	\begin{remark}
		The regularity theory for mixed operators  shows that \(w\) is a classical
		solution (see, e.g., \cite{Bia1,SVWZ} and the references therein for the case \(s_m=1\)), 
		we do not need this fact. For the regularity of solutions of
		equations involving the fractional Laplacian alone, see also \cite{Bernard}.
	\end{remark}

\section{Problem setting for the reduction and notation}\label{section.psn}

	From now on, \(s\in(0,1)\) and \(n\ge5\) with \(n>6-4s\) are fixed, and we write
	\[
		\tau\coloneqq\delta^{2(1-s)},\quad q_0\coloneqq\frac{2n}{n+2}=(2^*)',\quad
		f_\varepsilon(t)\coloneqq n(n-2)(t_+)^{p_\varepsilon},\quad p_\varepsilon=p_1-\varepsilon.
	\]
	Note that \(1<p_\varepsilon<p_1\), and that \(p_1\le2\) if and only if \(n\ge6\). 
	The case \(n=5\), where \(p_1=\frac73\), is covered by Lemma~\ref{lem:elem} below.

	\subsection{Functional framework}
		For \(\tau>0\), let
		\begin{align}\label{elespacio}
			\HH\coloneqq\Big\{u\in L^{2^*}(\R^n)\colon \nabla u\in L^2(\R^n),\
			[u]_s<\infty\Big\},\\
			\|u\|_{\tau,s}^2\coloneqq\|\nabla u\|_{L^2}^2+\tau[u]_s^2,
		\end{align}
		with inner product \(\langle u,v\rangle_{\tau,s}\coloneqq \langle u,v\rangle_{1}+\tau\langle u,v\rangle_s\).
		For every \(\tau>0\) this is a Hilbert space, and the norms \(\|\cdot\|_{\tau,s}\)
		are equivalent to each other, though not uniformly as \(\tau\to0\). We write
		\(\HH_{\rm rad}\) for the closed subspace of radial functions, \(\HH_{\rm rad}'\)
		for its dual, and \(\|\cdot\|_{\tau,s,*}\) for the dual norm. 
		
		Let \(g\in\HH'.\) A function \(u\in\HH\) is a weak solution of 
		\[ 
			-\Delta u+\tau(-\Delta)^su=g,
		\] 
		if \(\langle u,v\rangle_{\tau,s}=g(v)\) for every \(v\in\HH\).
		
		\begin{remark}\label{rem:symcrit}
			If \(u\in\HH_{\rm rad}\) and \(g\) is rotation invariant (for instance
			\(g=h\) or \(g=hu\) with \(h\) radial), it suffices to test with \(v\in\HH_{\rm rad}\):
			for \(v\in\HH\), its average \(\bar v\) over the orthogonal group belongs to
			\(\HH_{\rm rad}\) and \(\langle u,v\rangle_{\tau,s}=\langle u,\bar v\rangle_{\tau,s}\),
			\(g(v)=g(\bar v)\).
		\end{remark}
	
		By the Sobolev inequalities, for every \(u\in\HH\)
		\begin{equation}\label{eq:sob}
			\|u\|_{L^{2^*}}\le C_*\|\nabla u\|_{L^2}\le C_*\|u\|_{\tau,s},\qquad
			\|u\|_{L^{2^*_s}}\le C_{*,s}[u]_s\le C_{*,s}\tau^{-1/2}\|u\|_{\tau,s},
		\end{equation}
		where \(C_*=S^{-1/2}\) and \(C_{*,s}=S_s^{-1/2}.\) Here \(S\) and \(S_s\) are the best
		constants in 
		\[
			S\|u\|_{L^{2^*}}^2\le\|\nabla u\|_{L^2}^2 \quad \text{ and }\quad
			S_s\|u\|_{L^{2^*_s}}^2\le[u]_s^2.
		\]
		In particular, every \(h\in L^{q_0}(\R^n)\) defines an element of \(\HH'\) with
		\begin{equation}\label{eq:dualq0}
			\|h\|_{\tau,s,*}\le C_*\|h\|_{L^{q_0}},
		\end{equation}
		and for \(W\in L^{n/2}(\R^n)\) and \(\phi\in\HH\),
		\[
			\|W\phi\|_{\tau,s,*}\le C_*^2\|W\|_{L^{n/2}}\|\phi\|_{\tau,s}, 
		\]
		uniformly in \(\tau\).
		
		\medskip
		
		Notice that \(U\in\HH\), since \(U\in H^1(\R^n)\) for \(n\ge5\).
		
		\medskip
	
		We will use two auxiliary facts. The first one allows us to verify identities
		for test functions in \(\HH\cap L^2(\R^n)\), where the Fourier transform is available.
	
		\begin{lemma}\label{lem:trunc}
			For \(v\in\HH\) and \(k\in\mathbb{N}\), let
			\(T_kv\coloneqq\operatorname{sign}(v)\,(|v|-1/k)_+\). Then, for any \(\tau>0,\)
			\(T_kv\in\HH\cap L^2(\R^n)\) and \(T_kv\to v\) in \(\HH\)  as \(k\to\infty\).
		\end{lemma}
		\begin{proof}
			Observe that \(T_kv\) vanishes outside \(\{|v|>1/k\}\), a set of finite measure, and there
			\[
				|T_kv|^2\le|v|^2\le k^{2^*-2}|v|^{2^*}.
			\] 
			Hence \(T_kv\in L^2(\R^n)\). 
			
			On the other hand, since \(T_k\) is \(1\)-Lipschitz, 
			\[
				\|\nabla T_kv\|_{L^2}\le\|\nabla v\|_{L^2}\quad\text{ and  }\quad [T_kv]_s\le[v]_s. 
			\]
			Moreover
			\[v-T_kv=\max\{-1/k,\min\{v,1/k\}\}\eqqcolon c_k(v),\]
			so \(\nabla(v-T_kv)=\nabla v\,\chi_{\{|v|\le1/k\}}\to0\) in \(L^2(\R^n)\) (because \(\nabla v=0\) a.e.
			on \(\{v=0\}\)), and
			\[\left|c_k(v)(x)-c_k(v)(y)\right|\le\min\left\{|v(x)-v(y)|,\frac2k\right\},\]
			so \([v-T_kv]_s\to0\) by dominated
			convergence in the double integral. By \eqref{eq:sob}, also \(T_kv\to v\) in \(L^{2^*}(\R^n)\).
		\end{proof}
	
		The second fact gives the decay of \((-\Delta)^sU\), which is needed to obtain an
		error of order \(\tau\) (and not \(\tau^{1/2}\)) in the reduction.
	
		\begin{lemma}\label{lem:fracU}
			Let \(n\ge5\). The function \((-\Delta)^sU\), defined pointwise by
			\[
				(-\Delta)^sU(x)=\frac{C_{n,s}}{2}\int_{\R^n}\frac{2U(x)-U(x+y)-U(x-y)}{|y|^{n+2s}}\,dy,
			\]
			satisfies 
			\[
				|(-\Delta)^sU(x)|\le C(1+|x|)^{2-n-2s}\quad \forall x\in\R^n.
			\] 
			In particular, \((-\Delta)^sU\in L^{q_0}(\R^n)\) whenever \(n>6-4s\), which holds
			under our standing assumptions, and
			\begin{equation}\label{eq:fracU}
				\int_{\R^n}(-\Delta)^sU\,v\,dx=\langle U,v\rangle_s\qquad\text{for every }v\in\HH .
			\end{equation}
		\end{lemma}
		
		\begin{proof}
			Since \(U=-f_{n-2}\) with \(f_\beta(x)\coloneqq-(1+|x|^2)^{-\beta/2}\), the closed-form
			expression of \cite[formulas~(4.2) and~(6.5)]{FerrariVerbitsky}, applied with
			\(\alpha=2s\) and \(\beta=n-2\), gives
			\[
				(-\Delta)^sU(x)=C(s,n)\,{}_2F_1\!\left(\frac{n+2s}2,s+\frac{n-2}2,\frac n2,-|x|^2\right),
			\]
			where \(C(s,n)\) is a positive constant and \({}_2F_1\) is the Gaussian
			hypergeometric function. By the integral representation
			\cite[(7.1)]{FerrariVerbitsky} (valid here since
			\(\frac n2-\left(s+\frac{n-2}2\right)=1-s\in(0,1)\)), this expression is finite and
			continuous for every \(x\in\R^n\), which gives the bound for \(|x|\le2\). For
			\(|x|>2\), Pfaff's transformation \cite[Corollary~4.1, formula~(4.4)]{FerrariVerbitsky}
			gives the asymptotic behaviour
			\[
				|(-\Delta)^sU(x)|\sim C|x|^{2-n-2s}\quad\text{as }|x|\to\infty,
			\]
			which proves the stated bound, and shows \((-\Delta)^sU\in L^{q_0}(\R^n)\) if and
			only if \((n-2+2s)q_0>n\), that is, \(n>6-4s\).

			\medskip

			We now turn to \eqref{eq:fracU}. Since \(n\ge5\), \(U\in H^2(\R^n)\), and the
			pointwise \((-\Delta)^sU\) then coincides with
			\(\mathcal{F}^{-1}(|\xi|^{2s}\widehat U)\in L^2(\R^n)\), while
			\[
				\langle U,v\rangle_s=\int_{\R^n}
				|\xi|^{2s}\widehat U\,\overline{\widehat v}\,d\xi
				\quad \forall v\in H^s(\R^n),
			\]
			see \cite[Propositions 3.3 and 3.4]{DiNezzaPV}.

			This identification is stated there for Schwartz functions, but it extends
			to \(H^2(\R^n)\) because both sides are bounded operators from \(H^2(\R^n)\)
			to \(L^2(\R^n)\). For the integral formula, this follows from Minkowski's
			inequality together with
			\[\|2u-u(\cdot+y)-u(\cdot-y)\|_{L^2}\le\min\{|y|^2\|D^2u\|_{L^2},4\|u\|_{L^2}\}.\]

			By Plancherel, \eqref{eq:fracU} holds for
			\(v\in\HH\cap L^2(\R^n)\subset H^s(\R^n)\). Both sides of \eqref{eq:fracU} are
			continuous in \(v\) with respect to \(\|\cdot\|_{\tau,s}\): the left one by
			H\"older's inequality and \eqref{eq:sob}, the right one since
			\(|\langle U,v\rangle_s|\le[U]_s[v]_s\). Hence Lemma~\ref{lem:trunc} gives
			\eqref{eq:fracU} for every \(v\in\HH\).
		\end{proof}
	
		The following interpolation inequality is the point where the nonlocal term
		is used to control the nonlinearity at infinity.
	
		\begin{lemma}\label{lem:interp}
			Let \(r\in[2^*_s,2^*]\) and \(\theta(r)\in[0,1]\) be defined by
			\(\frac1r=\frac{\theta(r)}{2^*_s}+\frac{1-\theta(r)}{2^*}\). Then
			\[
				\|u\|_{L^r}\le C\,\tau^{-\theta(r)/2}\|u\|_{\tau,s}\qquad\text{for all }u\in\HH,\ \tau\in(0,1],
			\]
			with \(C=C(n,s)\). In particular, for \(r_\varepsilon\coloneqq p_\varepsilon q_0\)
			we have \(r_\varepsilon\in(2^*_s,2^*)\) for \(\varepsilon\) small,
			\[
				\theta_\varepsilon\coloneqq\theta(r_\varepsilon)=\frac{(n-2)\,\varepsilon}{2(1-s)\,p_\varepsilon},
			\]
			and, if \(\Lambda=[\lambda_-,\lambda_+]\subset(0,\infty)\) is a fixed compact interval,
			there are \(\varepsilon_1,C_\Lambda>0\) such that
			\begin{equation}\label{eq:key}
				\|u\|_{L^{r_\varepsilon}}\le C_\Lambda\|u\|_{\tau,s}\qquad
				\text{whenever } \tau=\lambda\varepsilon,\ \lambda\in\Lambda,\ \varepsilon\in(0,\varepsilon_1).
			\end{equation}
		\end{lemma}
		\begin{proof}
			The first inequality follows from H\"older's inequality and \eqref{eq:sob}.
			The formula for \(\theta_\varepsilon\) follows from
			\(\frac1{r_\varepsilon}-\frac1{2^*}=\frac{(n-2)\varepsilon}{2np_\varepsilon}\) and
			\(\frac1{2^*_s}-\frac1{2^*}=\frac{1-s}{n}\). Finally,
			\(\tau^{-\theta_\varepsilon/2}=\exp\left(\tfrac{\theta_\varepsilon}{2}|\log(\lambda\varepsilon)|\right)\)
			is bounded since \(\varepsilon|\log\varepsilon|\to0\).
		\end{proof}
	
		\begin{remark}
			Estimate \eqref{eq:key} explains why the reduction is carried out in the
			regime \(\tau\sim\varepsilon\). It also explains why no estimate in
			\(L^{q_0}(\R^n)\) for the solutions is needed (or available) uniformly in
			\(\tau\): the global integrability of the slightly subcritical nonlinearity
			is provided by the fractional part of the energy.
		\end{remark}

	\subsection{Positivity}

		We prove positivity in two steps: first \(u\ge0\), by testing the equation with
		the negative part, and then \(u>0\), via a self-contained argument based on the
		positivity of the Green function of \(-\Delta+\tau(-\Delta)^s\).

		\begin{lemma}\label{lem:positivity}
			If \(u\in\HH\) is a
			non-trivial weak solution of
			\[
				-\Delta u+\tau(-\Delta)^su=f_\varepsilon(u)
				\quad\text{ for some }\tau>0,
			\]
			then \(u>0\) a.e. in \(\R^n.\)
		\end{lemma}
		\begin{proof}
			\emph{Step 1: \(u\) is nonnegative.}
			Testing the equation with \(u_{-}\coloneqq \max\{-u,0\}\in\HH\) (where \(u=u_{+}-u_{-}\), \(u_{+}=\max\{u,0\}\)), and using \(f_\varepsilon(u)u_-=0\), we get
			\begin{equation}
				\label{weak}
				0=-\|\nabla u_-\|_{L^2}^2+\tau\langle u_{+}-u_{-},u_{-}\rangle_{s}
				=-\|\nabla u_-\|_{L^2}^2+\tau\left(\langle u_{+},u_{-}\rangle_{s} -[u_{-}]^{2}_{s}\right).
			\end{equation}
			Since 
			\[
				(u_+(x)-u_+(y))(u_-(x)-u_-(y))=-u_+(x)u_-(y)-u_+(y)u_-(x)\le0,
			\]
			we have
			\(\langle u_{+},u_{-}\rangle_{s}\le 0.\) Then, the three terms in \eqref{weak} are
			nonpositive and thus all vanish. Hence \(\nabla u_-=0\), and since \(u_-\in L^{2^*}(\R^n)\),
			\(u_-\equiv0\).
			
			\medskip
			
			\noindent\emph{Step 2: \(u\) is positive.} We split the argument into the following steps.

			\smallskip
			
			\noindent (i) \emph{Weak comparison.} 
			If \(g\in L^{q_0}(\R^n)\), \(g\ge0\), and \(w\in\HH\) is a
			weak solution of 
			\[
				-\Delta w+\tau(-\Delta)^sw=g,
			\] 
			then \(w\ge0\). 
			Indeed, testing with \(w_-\) and arguing as above,
			\[
				0\le\int_{\R^n} gw_- dx=-\|\nabla w_-\|_{L^2}^2+\tau\left(\langle w_+,w_-\rangle_s-[w_-]_s^2\right)\le-\|w_-\|_{\tau,s}^2.
			\]
			Hence \(\|w_-\|_{\tau,s}=0\), so \(w_-\equiv0\) and \(w\ge0\).

			\smallskip
			
			\noindent (ii) \emph{Positive solutions for nice data.} Let \(m(\xi)\coloneqq|\xi|^2+\tau|\xi|^{2s}\)
			and 
			\[
				\mathcal{P}_t\coloneqq\mathcal{F}^{-1}(e^{-tm})=\mathcal{G}_t*\mathcal{Q}_{\tau t}, \quad t>0, 
			\]
			where \(\mathcal{G}_t\) is the Gaussian kernel with 
			\(\widehat{\mathcal{G}_t}=e^{-t|\xi|^2}\) and \(\mathcal{Q}_{\tau t}\) 
			is the kernel with \(\widehat{\mathcal{Q}_{\tau t}}=e^{-\tau t|\xi|^{2s}}\). 
			By Bochner's subordination formula
			\[
				e^{-\lambda^s}=\int_0^\infty e^{-r\lambda}\nu_s(r)\,dr,\quad  \lambda\ge0, 
			\]
			with \(\nu_s\ge0\) and \(\int_0^\infty\nu_s\, dr=1\) (see \cite[Chapter IX.11]{Yosida}), 
			\(\mathcal{Q}_{\tau t}\) is a superposition of Gaussians, hence a probability density. 
			Since \(\mathcal{G}_t>0\) everywhere, 
			\[	
				\mathcal{P}_t(x)>0\quad\forall x, 
				\quad \int_{\R^n} \mathcal{P}_t dx=1\quad\text{ and }\quad 
				\|\mathcal{P}_t\|_\infty\le(4\pi t)^{-n/2}.
			\]

			Let \(0\le\tilde h\in L^1(\R^n)\cap L^\infty(\R^n)\), \(\tilde h\not\equiv0\), and
			\[
				\tilde w(x)\coloneqq\int_0^\infty(\mathcal{P}_t*\tilde h)(x)\,dt .
			\]

			Since \(0<\mathcal{P}_t*\tilde h\le\min\{\|\tilde h\|_\infty,(4\pi t)^{-n/2}\|\tilde h\|_1\}\) and \(n\ge3\),
			\(\tilde w\) is bounded and \(\tilde w(x)>0\) for every \(x\in\R^n\). For a Schwartz function \(\varphi\),
			\[
				\int_0^\infty e^{-tm}|\widehat{\tilde h}||\widehat\varphi|\,dt=
				\frac{|\widehat{\tilde h}||\widehat\varphi|}{m}\in L^1(\R^n),
			\]
			as \(|\widehat{\tilde h}|\le\|\tilde h\|_1\) and \(m^{-1}\le|\xi|^{-2}\) is locally integrable. Then,
			Fubini's theorem  gives \(\widehat{\tilde w}=\tfrac{\widehat{\tilde h}}{m}\) in the sense of distributions. 
			Hence
			\[
				\int_{\R^n} m|\widehat{\tilde w}|^2 \, d\xi=\int_{\R^n}\frac{|\widehat{\tilde h}|^2}{m} \, d\xi<\infty,
			\] 
			so \(\nabla\tilde w\in L^2(\R^n)\), and  \([\tilde w]_s<\infty\). Moreover \(\tilde w\in L^{2^*}(\R^n)\): splitting
			\[
				\widehat{\tilde w}=\widehat{\tilde w}\chi_{\{|\xi|\le1\}}+\widehat{\tilde w}\chi_{\{|\xi|>1\}},
			\] 
			the second piece is the transform of a function of \(H^1(\R^n)\subset L^{2^*}(\R^n)\). 
			The first is bounded by 
			\[ 
				\|\tilde h\|_1|\xi|^{-2}\chi_{\{|\xi|\le1\}}\in L^{\frac{2n}{n+2}}(\R^n),
			\]
			so by the Hausdorff--Young
			inequality it is the transform of a function of \(L^{2^*}(\R^n)\). Thus \(\tilde w\in\HH\).
			When \(n=3,4\), the function \(\tilde w\) need not belong to \(L^2(\R^n)\). In that case we first prove the identities below for \(\tilde w_\varepsilon\), defined by \(\widehat{\tilde w_\varepsilon}\coloneqq\widehat{\tilde h}/(m+\varepsilon)\), which belongs to \(H^1(\R^n)\). Since \(|\widehat{\tilde h}|/(m+\varepsilon)\le|\widehat{\tilde h}|/m\), dominated convergence allows us to let \(\varepsilon\to0\).
			In addition, for \(v\in\HH\cap L^2(\R^n)\), 
			\[
				\langle\tilde w,v\rangle_{\tau,s}=	\int_{\R^n} m\,\widehat{\tilde w}\,
				\overline{\widehat v}\, d\xi=\int_{\R^n}\widehat{\tilde h}\,\overline{\widehat v}\,d\xi=
				\int_{\R^n}\tilde hv \, dx,
			\]
			and by
			Lemma~\ref{lem:trunc} this holds for every \(v\in\HH\). So \(\tilde w\) is the weak
			solution with right-hand side \(\tilde h\).

			\smallskip

			\noindent \emph{Conclusion.} 
			Let \(h\coloneqq f_\varepsilon(u)=n(n-2)u^{p_\varepsilon}\ge0\), which is
		nontrivial. By Lemma~\ref{lem:interp}, \(u\in L^{r_\varepsilon}(\R^n)\), so 
		\(h\in L^{q_0}(\R^n)\). Let \(h_k\coloneqq\min\{h,k\}\chi_{B_k}\), which is nontrivial for \(k\) large, 
		and let \(w_k\) be given by (ii) with \(\tilde h=h_k\). Then \(u-w_k\in\HH\) is a weak solution with
		right-hand side \(h-h_k\in L^{q_0}(\R^n)\), \(h-h_k\ge0\), and by (i), \(u\ge w_k>0\) a.e.

		\smallskip
		Alternatively, the strict positivity can be obtained from a logarithmic lemma for
		mixed operators (see \cite{Ga}) arguing as in \cite[Lemma 3.3]{DelPezzoQuaas}.
	\end{proof}

	\subsection{The linearized operator}
	Let
	\[
		V\coloneqq f_0'(U)=n(n+2)U^{\frac4{n-2}}=\frac{n(n+2)}{(1+|x|^2)^2}\in L^{n/2}(\R^n).
	\]

		\begin{lemma}\label{lemaeq:ker}
			The space of radial solutions \(\phi\in\mathcal{D}^{1,2}(\R^n)\) of
			\(-\Delta\phi=V\phi\) is \(\spn\{\psi\}\),
			where
			\begin{equation}\label{la_psi}
							\psi(x)\coloneq x\cdot\nabla U+\frac{n-2}2 U=
							\frac{n-2}2\frac{1-|x|^2}{(1+|x|^2)^{n/2}}.
			\end{equation}
		\end{lemma}
		\begin{proof}
			See \cite{BianchiEgnell} or \cite{MR1696454}.
		\end{proof}
		
		\begin{remark}\label{remark.psi}
			Since \(\psi(x)=O(|x|^{2-n})\) and \(\nabla\psi(x)=O(|x|^{1-n})\), we have
			\(\psi\in H^1(\R^{n})\) if \(n\ge 5.\) 
			Hence \(\psi\in H^s(\R^n)\) and
			\(\psi\in\HH\), with \(\|\psi\|_{\tau,s}\le\|\psi\|_1\) for \(\tau\le1\).
			Moreover \(-\Delta\psi=V\psi\), so that
			\[
				\int_{\R^n}\nabla\phi\cdot\nabla\psi\,dx=\int_{\R^n} V\phi\psi\,dx\quad \forall \phi\in\HH.
			\]
			This is clear for \(\phi\in C^\infty_c(\R^n)\), both sides are continuous with
			respect to \(\|\nabla\phi\|_{L^2}\) (note that \(V\psi=O(|x|^{-n-2})\in L^{q_0}(\R^n)\)), 
			and \(C_c^\infty(\R^n)\) is dense in \(\mathcal{D}^{1,2}(\R^n)\).
		\end{remark}

		We define
		\[
			K\coloneqq\Big\{\phi\in\HH_{\rm rad}\colon
			\int_{\R^n}\nabla\phi\cdot\nabla\psi\,dx=0\Big\},
		\]
		which is a closed subspace of \(\HH_{\rm rad}\) for every \(\tau>0\).
		Given \(g\in\HH_{\rm rad}'\), we consider the projected linear problem
		\begin{equation}\label{eq:lin}
			\phi\in K,\ c\in\R\colon
			\langle\phi,v\rangle_{\tau,s}-\int_{\R^n}V\phi v\,dx=g(v)+c\int_{\R^n}V\psi v\,dx
			\quad\forall\,v\in\HH_{\rm rad}.
		\end{equation}

		\begin{lemma}\label{orz}
			There exist \(\tau_0>0\) and \(C>0\) such that for every \(\tau\in(0,\tau_0)\)
			and every \(g\in\HH_{\rm rad}'\), problem \eqref{eq:lin} has a unique solution
			\((\phi,c)\eqqcolon(\mathcal{A}_{\tau,s} g,c_{\tau,s}(g))\), and
			\[
				\|\phi\|_{\tau,s}+|c|\le C\|g\|_{\tau,s,*}.
			\]
		\end{lemma}
		\begin{proof}
			\emph{A priori estimate.} Assume by contradiction that there are
			\(\tau_k\to0\), \(g_k\) with \(\|g_k\|_{\tau_k,s,*}\to0\) and solutions
			\((\phi_k,c_k)\) of \eqref{eq:lin} with \(\|\phi_k\|_{\tau_k,s}+|c_k|=1\).
			Testing with \(v=\psi\) and using \(\phi_k\in K\) and Remark~\ref{remark.psi}
			(so that \(\int_{\R^n}\nabla\phi_k \cdot\nabla\psi \, dx=\int_{\R^n} 
			V\phi_k\psi \, dx=0\)),
			\[
				c_k\int_{\R^n}V\psi^2\,dx=\tau_k\langle\phi_k,\psi\rangle_s-g_k(\psi),
			\]
			and 
			\[
				|\tau_k\langle\phi_k,\psi\rangle_s|\le\tau_k^{1/2}\|\phi_k\|_{\tau_k,s}[\psi]_s\to0,
				\quad |g_k(\psi)|\le\|g_k\|_{\tau_k,s,*}\|\psi\|_1\to0. 
			\]
			Hence \(c_k\to0\), and since \(\|\phi_k\|_{\tau_k,s}+|c_k|=1\) for every \(k\),
			\(\|\phi_k\|_{\tau_k,s}\to1\).
			
			Up to a subsequence, \(\phi_k\rightharpoonup\phi\)
			weakly in \(\mathcal{D}^{1,2}(\R^n)\), with \(\phi\) radial and
			\[
				\int_{\R^n}\nabla\phi\cdot\nabla\psi dx=0.
			\] 
			For a fixed radial \(v\in C^\infty_c(\R^n),\) we have 
			\[
				\tau_k|\langle\phi_k,v\rangle_s|\le\tau_k^{1/2}[v]_s\to0,\quad
				g_k(v)\to0\quad\text{ and }\quad \int_{\R^n} V\phi_kv \, dx\to\int_{\R^n} V\phi v \,dx. 
			\]
			Therefore, tested with radial functions (and hence with all functions by averaging), we get 
			\[
				-\Delta\phi=V\phi,
			\]
			so that, by
			Lemma~\ref{lemaeq:ker} and the orthogonality, \(\phi=0\).

			Since \(V^{1/2}\in L^n(\R^n)\), multiplication by \(V^{1/2}\) is compact from
			\(\mathcal{D}^{1,2}(\R^n)\) into \(L^2(\R^n)\) (Rellich's theorem on balls and
			\(\|V^{1/2}\|_{L^n(\{|x|>R\})}\to0\)). Thus 
			\[
				\int_{\R^n} V\phi_k^2\to0.
			\]
			Testing
			\eqref{eq:lin} with \(v=\phi_k\) and using \(\int V\psi\phi_k=0\) we obtain
			\[
				\|\phi_k\|_{\tau_k,s}^2=\int_{\R^n}V\phi_k^2\,dx+g_k(\phi_k)\to0,
			\]
			a contradiction.
			
			\smallskip

			\noindent \emph{Existence.} Fix \(\tau\in(0,\tau_0)\). Let \(R_{\tau,s}\colon\HH_{\rm rad}'\to\HH_{\rm rad}\)
			be the Riesz isomorphism for \(\langle\cdot,\cdot\rangle_{\tau,s}\), \(\mathcal{B}\phi\coloneqq R_{\tau,s}(V\phi)\) and
			\(Z\coloneqq R_{\tau,s}(V\psi)\). By the compactness just mentioned and
			\[
				\|V w\|_{L^{q_0}}\le\|V^{1/2}\|_{L^n}\|V^{1/2}w\|_{L^2},
			\]
			\(\mathcal{B}\) is compact on \(\HH_{\rm rad}\).
			Since
			\[
				\langle\phi,Z\rangle_{\tau,s}=\int_{\R^n} V\psi\phi dx=\int_{\R^n}\nabla\psi\cdot\nabla\phi \, dx,
			\]
			we have \(K=Z^{\perp}\) (orthogonal complement in \(\langle\cdot,\cdot\rangle_{\tau,s}\)),
			and \eqref{eq:lin} is equivalent to \(\Pi_{\tau,s}(\phi-\mathcal{B}\phi)=\Pi_{\tau,s} R_{\tau,s} g\),
			\(\phi\in K\), where \(\Pi_{\tau,s}\) is the orthogonal projection onto \(K\)
			(the constant \(c\) is then the coefficient of \(Z\)). The operator
			\(\phi\mapsto\phi-\Pi_{\tau,s} \mathcal{B}\phi\) is a compact perturbation of the identity on \(K\),
			and it is injective by the a priori estimate. By the Fredholm alternative it is
			onto, which gives existence. 
			
			Uniqueness and the estimate follow from the a priori bound.
		\end{proof}

		\begin{remark}
			The constant in Lemma~\ref{orz} does not depend on \(\tau\), and the estimate
			involves only the energy norm \(\|\cdot\|_{\tau,s}\). In particular, the dilations
			\[
				u_R(x)\coloneqq R^{-n/q_0}\eta(x/R),
			\] 
			(which have constant \(L^{q_0}\)-norm but \(\|u_R\|_{\tau,s}\to0\)) are not an obstruction in this space.
		\end{remark}

	\subsection{Two integrals}

	We conclude this section with the following lemma, which determines the
	signs of the two quantities that drive the reduction.

	\begin{lemma} \label{taqui}
		Let \(n\ge 5\) and \(s\in(0,1)\). Then
		\begin{equation}\label{Pistoia}
			A\coloneqq n(n-2)\int_{\R^n} \log(U)\, U^{p_1} \psi \, dx
			=\frac{(n-2)^3}{4}\int_{\R^n}U^{2^*}dx>0,
		\end{equation}
		and
		\begin{equation}\label{chile241106}
			B\coloneqq\langle U,\psi \rangle_s=(s-1)[U]_s^2<0,
		\end{equation}
		where \(\psi\) was given in \eqref{la_psi}.
	\end{lemma}

	\begin{proof}
		Let \(U_t(x)\coloneqq t^{\frac{n-2}{2}}U(tx)\), \(t>0\), so that
		\(\frac{d}{dt}U_t\big|_{t=1}=\psi\). For \(n\ge5\) all the quantities below
		are finite and the differentiation under the integral sign is justified by
		the decay of \(U\) and \(\psi\).
		
		By scaling, 
		\[
			[U_t]_s^2=t^{2s-2}[U]_s^2.\] 
		Differentiating at \(t=1\) gives
		\[
			2\langle U,\psi\rangle_s=(2s-2)[U]_s^2, 
		\]
		which is \eqref{chile241106}. On the other hand, since for \(q\) in a neighborhood of \(2^*\), 
		\[	\int_{\R^n}U_t^q\,dx=t^{\frac{(n-2)q}{2}-n}\int_{\R^n}U^q\,dx,\]
		differentiating again at \(t=1\) it follows that
		\begin{equation}\label{IA}
			q\int_{\R^n}U^{q-1}\psi\,dx=\Big(\frac{(n-2)q}{2}-n\Big)\int_{\R^n}U^q\,dx,
		\end{equation}
		which, by choosing \(q=2^*\), clearly implies 
		\[
			\int_{\R^n} U^{p_1}\psi \, dx=0.
		\]
		Differentiating \eqref{IA} with respect to \(q\) and evaluating later at \(q=2^*\) the previous identity yields
		\[
			2^*\int_{\R^n}\log(U)U^{p_1}\psi \, dx=\frac{n-2}{2}\int_{\R^n} U^{2^*}\, dx,
		\]
		which gives \eqref{Pistoia}.
	\end{proof}

	In particular, the number \(\lambda_0\) in Theorem~\ref{ElTeorema} is
	\[
		\lambda_0=-\frac{A}B=\frac{(n-2)^3}{4(1-s)}\dfrac{\displaystyle \int_{\R^n} U^{2^*}\, dx}{[U]_s^2}.
	\]

\section{Existence}\label{sec:existence}

	Throughout this section we fix a compact interval
	\(\Lambda=[\lambda_-,\lambda_+]\) with \(0<\lambda_-<\lambda_0<\lambda_+\) and take
	\[
		\tau=\delta^{2(1-s)}=\lambda\varepsilon,\qquad\lambda\in\Lambda .
	\]
	All constants below are independent of \(\varepsilon\) small and \(\lambda\in\Lambda\).

	\subsection{The nonlinear projected problem}
	We look for \(\phi\in K\) and \(c\in\R\) such that
	\begin{equation}\label{eq:redu}
		-\Delta(U+\phi)+\tau(-\Delta)^s(U+\phi)=f_\varepsilon(U+\phi)+c\,V\psi
	\end{equation}
	weakly in \(\HH_{\rm rad}\). Using \(-\Delta U=f_0(U)\), \eqref{eq:redu} can be written as
	\begin{equation}\label{eq:redu2}
		\langle\phi,v\rangle_{\tau,s}-\int_{\R^n} V\phi v\,dx=\mathcal{E}(v)+\int_{\R^n} N_\varepsilon(\phi)v\,dx
		+c\int_{\R^n} V\psi v\,dx
	\end{equation}
	for all \( v\in\HH_{\rm rad},\) where
	\[
		\mathcal{E}(v)\coloneqq\int_{\R^n}\left(f_\varepsilon(U)-f_0(U)\right)v\,dx-\tau\langle U,v\rangle_s,
	\]
	\[
		N_\varepsilon(\phi)\coloneqq\left[f_\varepsilon(U+\phi)-f_\varepsilon(U)-f_\varepsilon'(U)\phi\right]
		+\left[f_\varepsilon'(U)-f_0'(U)\right]\phi .
	\]

	\begin{lemma}\label{lem:elem}
		Let \(1<p_-<p_+<\infty\). There is \(C=C(p_-,p_+)\) such that for \(a>0\), \(p\in[p_-,p_+]\)
		and \(\Theta(b)\coloneqq(a+b)_+^p-a^p-pa^{p-1}b\), with \(\chi_p\coloneqq1\) if \(p>2\) and
		\(\chi_p\coloneqq0\) if \(p\le2\),
		\begin{equation}\label{eq:elem}
		\begin{gathered}
			|\Theta(b)|\le C\left(|b|^p+\chi_pa^{p-2}b^2\right),\\
			|\Theta(b_1)-\Theta(b_2)|\le C\left(|b_1|^{p-1}+|b_2|^{p-1}+\chi_pa^{p-2}(|b_1|+|b_2|)\right)|b_1-b_2| .
		\end{gathered}
		\end{equation}
	\end{lemma}
	\begin{proof}
		We have \(\Theta(0)=0\) and \(\Theta'(t)=p\left((a+t)_+^{p-1}-a^{p-1}\right)\). 
		
		If \(p\le2\), \(x\mapsto x_+^{p-1}\) is \((p-1)\)-H\"older continuous with constant \(1\), so
		\[
			|\Theta'(t)|\le p|t|^{p-1}.
		\] 
		
		If \(p>2,\) the Mean Value Theorem gives 
		\[
			|\Theta'(t)|\le
			\begin{cases}
				p(p-1)\left(\frac{3a}2\right)^{p-2}|t| &\text{if } |t|\le \frac{a}2,\\
				p\left((a+|t|)^{p-1}+a^{p-1}\right)\le C|t|^{p-1} &\text{if } |t|>\frac{a}2.
			\end{cases}
		\] 
		Hence
		\[	
			|\Theta'(t)|\le C(|t|^{p-1}+\chi_pa^{p-2}|t|),
		\]
		and \eqref{eq:elem} follows by integration.
	\end{proof}

	\begin{lemma}\label{lemma:estimates}
		For \(\varepsilon\) small, \(\lambda\in\Lambda\) and \(\phi,\phi_1,\phi_2\in\HH\),
		\begin{enumerate}[label=\textup{(\alph*)}]
			\item \(\|\mathcal{E}\|_{\tau,s,*}\le C(\varepsilon+\tau)\le C\varepsilon\),
			\item \(\|N_\varepsilon(\phi)\|_{\tau,s,*}\le C\left(\|\phi\|_{\tau,s}^{p_\varepsilon}+\|\phi\|_{\tau,s}^2+\varepsilon\|\phi\|_{\tau,s}\right)\),
			\item \(\begin{aligned}[t]
				\|N_\varepsilon(\phi_1)&-N_\varepsilon(\phi_2)\|_{\tau,s,*}\le\\
				&C\left(\|\phi_1\|_{\tau,s}^{p_\varepsilon-1}+\|\phi_2\|_{\tau,s}^{p_\varepsilon-1}+\|\phi_1\|_{\tau,s}+\|\phi_2\|_{\tau,s}+\varepsilon\right)\|\phi_1-\phi_2\|_{\tau,s}.
			\end{aligned}\)
		\end{enumerate}
	\end{lemma}
	\begin{proof}
		Recall that \(0<U\le1\), \(U(x)\asymp(1+|x|)^{2-n}\) and
		\[
			|\log U|=\frac{n-2}{2}\log(1+|x|^2)\le C\log(2+|x|).
		\]

		\noindent (a) Since 
		\[
			U^{-\varepsilon}-1=e^{\varepsilon|\log U|}-1\le\varepsilon|\log U|e^{\varepsilon|\log U|},
		\] 
		we get
		\[
			|U^{p_1-\varepsilon}-U^{p_1}|\le\varepsilon|\log U|U^{p_1-\varepsilon}
			\le C\varepsilon\log(2+|x|)(1+|x|)^{-(n+2)+(n-2)\varepsilon},
		\] 
		whose \(L^{q_0}\) norm is \(O(\varepsilon)\) since \((n+2)q_0=2n\). 
		By \eqref{eq:dualq0}, the first term of \(\mathcal{E}\)
		has dual norm \(O(\varepsilon)\). For the second one,
		by Lemma~\ref{lem:fracU},
		\[
			\tau|\langle U,v\rangle_s|=\tau\left|\int_{\R^n}(-\Delta)^sU\,v\,dx\right|\le
			\tau\|(-\Delta)^sU\|_{L^{q_0}}C_*\|v\|_{\tau,s}.
		\]

		\noindent(b) We use Lemma~\ref{lem:elem} above with \(a=U\), \(p=p_\varepsilon\). The first
		bracket in \(N_\varepsilon(\phi)\) is bounded by \(C(|\phi|^{p_\varepsilon}+U^{p_\varepsilon-2}\phi^2)\), the
		second term being present only if \(p_\varepsilon>2\), i.e. \(n=5\). For it we use
		H\"older's inequality with exponents \((\frac{2n}{6-n},\frac{2^*}{2},2^*)\) and
		\(U^{p_\varepsilon-2}\in L^{\frac{2n}{6-n}}(\R^n)\) (for \(n=5\), \((n-2)(p_\varepsilon-2)\frac{2n}{6-n}\to10>5\)),
		which gives 
		\[
			\int_{\R^n} U^{p_\varepsilon-2}\phi^2|v|\,dx\le C\|\phi\|_{\tau,s}^2\|v\|_{\tau,s}.
		\]
		For the first term,
		by H\"older's inequality (\(p_\varepsilon/r_\varepsilon+1/2^*=1\)), \eqref{eq:key} and \eqref{eq:sob},
		\[
			\int_{\R^n}|\phi|^{p_\varepsilon}|v|\,dx\le\|\phi\|_{L^{r_\varepsilon}}^{p_\varepsilon}\|v\|_{L^{2^*}}
			\le C\|\phi\|_{\tau,s}^{p_\varepsilon}\|v\|_{\tau,s} .
		\]
		For the second bracket, writing
		\[
			f_\varepsilon'(U)-f_0'(U)=n(n-2)\left[(p_\varepsilon-p_1)U^{p_\varepsilon-1}+p_1(U^{p_\varepsilon-1}-U^{p_1-1})\right],
		\]
		and arguing as in (a), we get
		\begin{align*}
			W_\varepsilon&\coloneqq|f_\varepsilon'(U)-f_0'(U)|\\
			&\le
			C\varepsilon(1+|\log U|)U^{p_1-1-\varepsilon}\le C\varepsilon\log(2+|x|)(1+|x|)^{-4+(n-2)\varepsilon}.
		\end{align*}
		Thus 
		\[
			\|W_\varepsilon\|_{L^{n/2}}\le C\varepsilon\quad\text{ and }
			\|W_\varepsilon\phi\|_{\tau,s,*}\le C\varepsilon\|\phi\|_{\tau,s}.
		\]

		\noindent(c) We use the second inequality in \eqref{eq:elem}, the bound for \(W_\varepsilon\), the
			previous H\"older inequality for the term \(U^{p_\varepsilon-2}(|\phi_1|+|\phi_2|)|\phi_1-\phi_2|\), and H\"older's inequality
		with \(\frac{p_\varepsilon-1}{r_\varepsilon}+\frac1{r_\varepsilon}+\frac1{2^*}=1\) together with
		\eqref{eq:key}.
	\end{proof}

	\begin{lemma}\label{eq:puntofijo}
		There exist \(\varepsilon_0,C_0>0\) such that for \(\varepsilon\in(0,\varepsilon_0)\) and
		\(\lambda\in\Lambda\), problem \eqref{eq:redu} has a unique solution
		\((\phi_{\varepsilon,\lambda},c_{\varepsilon,\lambda})\) with \(\phi_{\varepsilon,\lambda}\in K\) and
		\(\|\phi_{\varepsilon,\lambda}\|_{\tau,s}\le C_0\varepsilon\). Moreover, for fixed \(\varepsilon\), the maps
		\(\lambda\mapsto\phi_{\varepsilon,\lambda}\in\HH\) and \(\lambda\mapsto c_{\varepsilon,\lambda}\) are continuous.
	\end{lemma}
	\begin{proof}
		By Lemma~\ref{orz} (note \(\tau\le\lambda_+\varepsilon<\tau_0\)), \eqref{eq:redu2} is
		equivalent to the fixed point problem
		\[
			\phi=T(\phi)\coloneqq\mathcal{A}_{\tau,s}\left(\mathcal{E}+N_\varepsilon(\phi)\right) \text{ in } K,
		\] 
		and then
		\(c=c_{\tau,s}(\mathcal{E}+N_\varepsilon(\phi))\). Let
		\(X_\varepsilon\coloneqq\{\phi\in K\colon\|\phi\|_{\tau,s}\le C_0\varepsilon\}\). By Lemmas~\ref{orz} and
		\ref{lemma:estimates},
		\[
			\|T(\phi)\|_{\tau,s}\le C_1\left(\varepsilon+(C_0\varepsilon)^{p_\varepsilon}+(C_0\varepsilon)^2+C_0\varepsilon^2\right)\le C_0\varepsilon,
		\]
		if \(C_0=2C_1\) and \(\varepsilon\) is small, because \(p_\varepsilon\ge\frac{1+p_1}{2}>1\). Likewise,
		\[
			\|T(\phi_1)-T(\phi_2)\|_{\tau,s}\le C_1\left(2(C_0\varepsilon)^{p_\varepsilon-1}+2C_0\varepsilon+\varepsilon\right)\|\phi_1-\phi_2\|_{\tau,s}
			\le\tfrac12\|\phi_1-\phi_2\|_{\tau,s},
		\]
		on \(X_\varepsilon\). The contraction mapping principle gives existence and uniqueness.

		For fixed \(\varepsilon\) the norms \(\|\cdot\|_{\lambda\varepsilon}\), \(\lambda\in\Lambda\), are
		uniformly equivalent, and the estimates above hold with the same constants on
		the larger balls \(\{\|\phi\|_{\tau,s}\le2C_0\varepsilon\}\), which contain \(\phi_{\lambda'}\) for
		\(\lambda'\) close to \(\lambda\). If \(\lambda_k\to\lambda\), then, writing \(T_\lambda\) for the map
		above, 
		\[
			\|\phi_{\lambda_k}-\phi_\lambda\|\le2\|T_{\lambda_k}\phi_\lambda-T_\lambda\phi_\lambda\|,
		\] 
		by the
		uniform contraction property. We affirm now that the right-hand side tends to \(0\) because
		\(\mathcal{E}\) depends affinely on \(\tau\) and \(\tau\mapsto\mathcal{A}_{\tau,s}g\) is
		continuous for fixed \(g\in\HH_{\rm rad}'\), uniformly for \(\tau\) in compact subsets of
		\((0,\tau_0)\). Indeed, if we denoted by \(\phi\coloneqq\mathcal{A}_{\tau,s}g\) and
		\(\phi'\coloneqq\mathcal{A}_{\tau',s}g\), subtracting the two instances of \eqref{eq:lin}
		shows that \(\phi'-\phi\in K\) solves \eqref{eq:lin} at parameter \(\tau\) with right-hand
		side \(v\mapsto-(\tau'-\tau)\langle\phi',v\rangle_s\), so, using
		\(|\langle\phi',v\rangle_s|\le[\phi']_s[v]_s\le\tau^{-1/2}[\phi']_s\|v\|_{\tau,s}\) and the a
		priori bound \([\phi']_s\le\tau'^{-1/2}\|\phi'\|_{\tau',s}\le C\tau'^{-1/2}\|g\|_{\tau',s,*}\),
		Lemma~\ref{orz} gives
		\[
			\|\phi'-\phi\|_{\tau,s}\le C\tau^{-1/2}|\tau'-\tau|\,[\phi']_s
			\le C(\tau\tau')^{-1/2}|\tau'-\tau|\,\|g\|_{\tau',s,*}.
		\]
		This gives the claimed continuity, with a constant uniform for \(\tau,\tau'\) in a fixed
		compact subinterval of \((0,\tau_0)\). The continuity of \(c_{\varepsilon,\lambda}\) in \(\lambda\) follows
		from the exact identity
		\[
			c\int_{\R^n}V\psi^2\,dx=\int_{\R^n}(f_0(U)-f_\varepsilon(U))\psi\,dx-\int_{\R^n}N_\varepsilon(\phi)\psi\,dx+\tau B+\tau\langle\phi,\psi\rangle_s,
		\]
		obtained in the proof of Lemma~\ref{lemaimp} below by testing \eqref{eq:redu} with \(\psi\), whose right-hand side
		depends continuously on \(\lambda\).
	\end{proof}

	\subsection{The reduced equation}
	We derive now the equation satisfied by the Lagrange multiplier
	\(c_{\varepsilon,\lambda}\), whose sign for \(\varepsilon\) small determines, via the
	Intermediate Value Theorem, the parameter \(\lambda_\varepsilon\) in Theorem~\ref{ElTeorema}.
	\begin{lemma} \label{lemaimp}
		As \(\varepsilon\to0\), uniformly in \(\lambda\in\Lambda\),
		\begin{equation}\label{eq:c}
			c_{\varepsilon,\lambda}\int_{\R^n}V\psi^2\,dx=\varepsilon\left(A+\lambda B+\omega(\varepsilon,\lambda)\right),
			\qquad \sup_{\lambda\in\Lambda}|\omega(\varepsilon,\lambda)|\to0,
		\end{equation}
		where \(A>0\) and \(B<0\) are given in Lemma~\ref{taqui}.
	\end{lemma}
	\begin{proof}
		We denote by \(\phi\coloneqq\phi_{\varepsilon,\lambda}\), \(c\coloneqq c_{\varepsilon,\lambda}\). Testing \eqref{eq:redu} with
		\(v=\psi\in\HH_{\rm rad}\), and using 
		\[
			\int_{\R^n}\nabla U\cdot\nabla\psi\,dx=\int_{\R^n} f_0(U)\psi\,dx \quad\text{ and }
			\int_{\R^n}\nabla\phi\cdot\nabla\psi\,dx=\int_{\R^n} f_0'(U)\phi\psi\,dx=0,
		\]
		since \(\phi\in K\) and
		Remark~\ref{remark.psi}, we obtain
		\begin{align*}
			c\int_{\R^n} V\psi^2\,dx&=-\int_{\R^n}\left[f_\varepsilon(U+\phi)-f_0(U)-f_0'(U)\phi\right]\psi\,dx
			+\tau\langle U,\psi\rangle_s+\tau\langle\phi,\psi\rangle_s\\
			&=\int_{\R^n}\left(f_0(U)-f_\varepsilon(U)\right)\psi\,dx-\int_{\R^n} N_\varepsilon(\phi)\psi\,dx
			+\tau B+\tau\langle\phi,\psi\rangle_s .
		\end{align*}
		Note that \(N_\varepsilon(\phi)\) contains the term \((f_\varepsilon'(U)-f_0'(U))\phi\). We estimate
		each adding term in the previous expression separately. First of all we observe that
			\[f_0(U)-f_\varepsilon(U)=n(n-2)U^{p_1}(1-U^{-\varepsilon}),\] and
			\[
				|1-U^{-\varepsilon}-\varepsilon\log U|\le\varepsilon^2(\log U)^2U^{-\varepsilon},\quad
				\int_{\R^n}(\log U)^2U^{p_1-\varepsilon}|\psi|\,dx\le C.
			\]
			So, recalling the definition of \(A\) in Lemma~\ref{taqui},
			\[
				\int_{\R^n}\left(f_0(U)-f_\varepsilon(U)\right)\psi\,dx
				=\varepsilon\,n(n-2)\int_{\R^n}U^{p_1}\log U\,\psi\,dx+O(\varepsilon^2)
				=\varepsilon A+O(\varepsilon^2).
			\]
			Secondly, by Lemma~\ref{lemma:estimates}(b) and Lemma~\ref{eq:puntofijo},
			\[
				\Big|\int_{\R^n} N_\varepsilon(\phi)\psi\,dx\Big|\le C\left((C_0\varepsilon)^{p_\varepsilon}+(C_0\varepsilon)^2+C_0\varepsilon^2\right)\|\psi\|_1=o(\varepsilon).
			\]
			Finally we get that $\tau B=\varepsilon\lambda B$ and $|\tau\langle\phi,\psi\rangle_s|\le\tau^{1/2}\|\phi\|_{\tau,s}[\psi]_s\le C\varepsilon^{3/2}$. These four estimates allow us to conclude \eqref{eq:c}.
	\end{proof}

	We can finally prove our main result.
	\begin{proof}[Proof of Theorem \ref{ElTeorema}]

		Since \(A>0>B\) and \(\lambda_-<\lambda_0=-A/B<\lambda_+\), we have
		\(A+\lambda_-B>0>A+\lambda_+B\). By \eqref{eq:c}, for \(\varepsilon\) small,
		\(c_{\varepsilon,\lambda_-}>0>c_{\varepsilon,\lambda_+}\), and by the continuity in
		Lemma~\ref{eq:puntofijo} there is \(\lambda_\varepsilon\in(\lambda_-,\lambda_+)\) with
		\(c_{\varepsilon,\lambda_\varepsilon}=0\). Moreover, \(c_{\varepsilon,\lambda_\varepsilon}=0\) and \eqref{eq:c} give
		\[
			|A+\lambda_\varepsilon B|\le\sup_{\lambda\in\Lambda}|\omega(\varepsilon,\lambda)|\to0,
		\] 
		so \(\lambda_\varepsilon\to-A/B=\lambda_0\).

		Let \(\delta_\varepsilon\coloneqq(\lambda_\varepsilon\varepsilon)^{1/(2(1-s))}\) and
		\(z_\varepsilon\coloneqq U+\phi_{\varepsilon,\lambda_\varepsilon}\). Then \(z_\varepsilon\) is a weak solution of
		\(-\Delta z+\delta_\varepsilon^{2(1-s)}(-\Delta)^sz=f_\varepsilon(z)\) tested with radial functions,
		hence with all \(v\in\HH\) by Remark~\ref{rem:symcrit}. It is nontrivial, so it is
		positive by Lemma~\ref{lem:positivity}, and
		\[
			\|\nabla\phi_\varepsilon\|_{L^2}+\delta_\varepsilon^{1-s}[\phi_\varepsilon]_s\le2\|\phi_\varepsilon\|_{\tau,s}\le2C_0\varepsilon.
		\]

		Finally, \(u_\varepsilon(y)=\delta_\varepsilon^{-2/(p_\varepsilon-1)}z_\varepsilon(y/\delta_\varepsilon)\) solves
		\eqref{elproblema} by the scaling described in the Introduction. Since
		\(\frac{4}{p_\varepsilon-1}=n-2+O(\varepsilon)\), a change of variables gives, for every \(\rho\ge0\),
		\[
			\int_{\{|y|>\rho\}}|\nabla u_\varepsilon|^2\,dy=\delta_\varepsilon^{\,n-2-\frac{4}{p_\varepsilon-1}}
				\int_{\{|x|>\rho/\delta_\varepsilon\}}|\nabla z_\varepsilon|^2\,dx,
		\]
		and \(\delta_\varepsilon^{O(\varepsilon)}\to1\) because \(\varepsilon|\log\delta_\varepsilon|\to0\). Since
		\(z_\varepsilon\to U\) in \(\mathcal{D}^{1,2}(\R^n)\), the right-hand side tends to
		\(\int_{\R^n}|\nabla U|^2\,dx\) for \(\rho=0\) and to \(0\) for every \(\rho>0\), which gives the
		concentration statement.
	\end{proof}

\section{Proof of Theorem~\ref{TeoremaBlowDown}}\label{sec:blowdown}

	Throughout this section \(s\in(0,1)\), \(n\ge3\), and
	\[
		q_s\coloneqq(2^*_s)'=\frac{2n}{n+2s},\qquad
		F_\eta(t)\coloneqq c_{n,s}(t_+)^{p_s+\eta},\]
	\[
		V_s\coloneqq F_0'(U_s)=c_{n,s}p_s\,U_s^{p_s-1}=\frac{c_{n,s}p_s}{(1+|x|^2)^{2s}} .
	\]
	Since 
	\[
		U_s\asymp(1+|x|)^{2s-n},
		\,|\nabla U_s|\le C(1+|x|)^{2s-n-1}\text{ and }
		|D^2U_s|\le C(1+|x|)^{2s-n-2},
	\]
	a direct
	computation shows that, for \(n\ge3\),
	\begin{equation}\label{eq:Us-int}
	\begin{gathered}
		U_s\in L^{2^*}(\R^n)\cap L^{2^*_s}(\R^n),\quad \nabla U_s\in L^2(\R^n),\\
		\Delta U_s\in L^{q_s}(\R^n)\cap L^{(2^*)'}(\R^n),\quad V_s\in L^{\frac{n}{2s}}(\R^n).
	\end{gathered}
	\end{equation}
	We clarify that, for instance, \(\Delta U_s\in L^{q_s}(\R^n)\) 
	because \((n+2-2s)q_s>n\), i.e. \(n>6s-4\), which in particular implies \(U_s\in\HH\).

	\subsection{Functional setting}
	We use the space \(\HH\) of \eqref{elespacio} (as a set), now with the norm and
	inner product
	\[
		\|u\|_{\sigma,s}^2\coloneqq[u]_s^2+\sigma\|\nabla u\|_{L^2}^2,\]
	\[
		\langle u,v\rangle_{\sigma,s}\coloneqq\langle u,v\rangle_s+\sigma\int_{\R^n}\nabla u\cdot\nabla v\,dx,
		\qquad\sigma\in(0,1].
	\]
	Since \(\|u\|_{\sigma,s}^2=\sigma\|u\|_{1/\sigma,s}^2\), \((\HH,\|\cdot\|_{\sigma,s})\) is a Hilbert space, and all
	the properties of \(\HH\) established in Section~\ref{section.psn}
	(Remark~\ref{rem:symcrit}, Lemma~\ref{lem:trunc}) remain valid. We denote by
	\(\|\cdot\|_{\sigma,s,*}\) the dual norm. By the Sobolev inequalities,
	\begin{equation}\label{eq:sob-s}
		\|u\|_{L^{2^*_s}}\le C[u]_s\le C\|u\|_{\sigma,s},\qquad
		\|u\|_{L^{2^*}}\le C\|\nabla u\|_{L^2}\le C\sigma^{-1/2}\|u\|_{\sigma,s}.
	\end{equation}
	Hence
	\[
		\|h\|_{\sigma,s,*}\le C\|h\|_{L^{q_s}}\quad\text{ and }\quad
		\|W\phi\|_{\sigma,s,*}\le C\|W\|_{L^{n/(2s)}}\|\phi\|_{\sigma,s},
	\] 
	uniformly in \(\sigma\).

	\begin{lemma}\label{lem:interp-s}
		Let \(r\in[2^*_s,2^*]\) and \(\frac1r=\frac{1-\theta}{2^*_s}+\frac{\theta}{2^*}\). Then
		\(\|u\|_{L^r}\le C\sigma^{-\theta/2}\|u\|_{\sigma,s}\) for \(u\in\HH\) and \(\sigma\in(0,1]\). For
		\(r_\eta\coloneqq(p_s+\eta)q_s\) we have
		\(\theta_\eta=\frac{(n-2s)\eta}{2(1-s)(p_s+\eta)}\), and if \(\Lambda=[\lambda_-,\lambda_+]\subset(0,\infty)\)
		then \(\|u\|_{L^{r_\eta}}\le C_\Lambda \|u\|_{\sigma,s}\) whenever \(\sigma=\lambda\eta\), \(\lambda\in\Lambda\),
		\(\eta\) small.
	\end{lemma}
	\begin{proof}
		H\"older's inequality and \eqref{eq:sob-s}, the value of \(\theta_\eta\) follows from
		\(\frac1{2^*_s}-\frac1{r_\eta}=\frac{(n-2s)\eta}{2n(p_s+\eta)}\) and
		\(\frac1{2^*_s}-\frac1{2^*}=\frac{1-s}{n}\), and \(\sigma^{-\theta_\eta/2}\) is bounded since
		\(\eta|\log\eta|\to0\). We note that
		Remark~\ref{rem:symcrit}, Lemmas~\ref{lem:trunc} and~\ref{lem:positivity} hold for every
		\(n\ge3\), and they are used in this section without change.
	\end{proof}

	\subsection{The linearized operator and two integrals}
	Let
	\[
		\psi_s(x)\coloneqq\frac{d}{dt}\Big(t^{\frac{n-2s}{2}}U_s(tx)\Big)\Big|_{t=1}
		=x\cdot\nabla U_s+\frac{n-2s}{2}U_s=\frac{n-2s}{2}\,\frac{1-|x|^2}{(1+|x|^2)^{\frac{n-2s+2}{2}}}.	\]
	Then $\psi_s\asymp|x|^{2s-n}$ and $|\nabla\psi_s|\le C|x|^{2s-n-1}$
	at infinity, so \(\psi_s\in\HH\) and 
	\[
		\|\psi_s\|_{\sigma,s}\le[\psi_s]_s+\|\nabla\psi_s\|_{L^2}. 
	\]
	Differentiating the equation
	satisfied by the dilations of \(U_s\) we get 
	\[
		(-\Delta)^s\psi_s=V_s\psi_s,
	\] 
	so that
	\[
		\langle\phi,\psi_s\rangle_s=\int_{\R^n} V_s\psi_s\phi\,dx\quad\forall\phi\in\HH.
	\]
	By the nondegeneracy of
	\(U_s\) (see \cite{DdPS}),  we know that the bounded radial solutions of 
	\[
		(-\Delta)^s\phi=V_s\phi,
	\] 
	are the
	multiples of \(\psi_s\). This nondegeneracy statement is phrased for
	\emph{bounded} radial solutions, but what we actually need is nondegeneracy for every
	radial solution \(\phi\in\dot H^s(\R^n)\), without assuming boundedness in advance. The
	extension is a standard bootstrap argument, which we spell out here. First, since
	\(V_s\in L^{n/(2s)}(\R^n)\cap L^\infty(\R^n)\), a Brezis--Kato type iteration (see
	\cite{DMV} for the same argument for the fractional Laplacian in \(\R^n\)) improves the
	integrability of \(\phi\) step by step, from \(\phi\in L^{2^*_s}(\R^n)\) (all that
	\(\dot H^s(\R^n)\) gives directly) up to \(\phi\in L^r(\R^n)\) for every
	\(r\in[2^*_s,\infty)\). Second, once this higher integrability is available we can
	represent \(\phi\) through the equation itself: writing
	\(\phi=(-\Delta)^{-s}(V_s\phi)=I_{2s}*(V_s\phi)\), where \(I_{2s}\) is the Riesz
	potential (see \cite{Stein}), it only remains to check that this convolution is
	bounded. This holds because \(V_s\phi\in L^{r_1}(\R^n)\cap L^{r_2}(\R^n)\) for some
	\(r_1<\frac n{2s}<r_2,\) having the right-hand side simultaneously in a space below and
	a space above the scaling exponent \(n/(2s)\) of \(I_{2s}\) is exactly what makes the
	convolution bounded, one exponent controlling the singularity of \(I_{2s}\) near the
	origin and the other controlling its slow decay at infinity. Hence \(\phi\) is bounded,
	and the nondegeneracy statement applies to it.

	\begin{lemma}\label{taqui-s}
		We have
		\[
			A_s\coloneqq c_{n,s}\int_{\R^n}\log(U_s)\,U_s^{p_s}\psi_s\,dx
			=c_{n,s}\frac{(n-2s)^2}{4n}\int_{\R^n}U_s^{2^*_s}\,dx>0,
		\]
		\[
			B_s\coloneqq\int_{\R^n}\nabla U_s\cdot\nabla\psi_s\,dx=\int_{\R^n}(-\Delta U_s)\psi_s\,dx
			=(1-s)\|\nabla U_s\|_{L^2}^2>0 .
		\]
	\end{lemma}
	\begin{proof}
		Let \((U_s)_t(x)\coloneqq t^{\frac{n-2s}{2}}U_s(tx)\), \(t>0\), so that
		\(\frac{d}{dt}(U_s)_t\big|_{t=1}=\psi_s\). As in Lemma~\ref{taqui}, all the quantities
		below are finite and the differentiation under the integral sign is justified by
		the decay of \(U_s\) and \(\psi_s\). By scaling,
		\[
			\|\nabla(U_s)_t\|_{L^2}^2=t^{2-2s}\|\nabla U_s\|_{L^2}^2 .
		\]
		Following the same idea as in Lemma~\ref{taqui} we derive at \(t=1\) to obtain \(B_s\). The corresponding integration by parts
		follows from \eqref{eq:Us-int} and the density of \(C^\infty_c\) in
		\(\mathcal{D}^{1,2}\). Indeed, the identity 
		\[
			\int_{\R^n}\nabla U_s\cdot\nabla v\,dx=\int_{\R^n}(-\Delta U_s)v\,dx
		\]
		holds for all \(v\in\HH\) by density, since \(|\Delta U_s|\le C(1+|x|)^{2s-n-2}\), so
		\[
			\Delta U_s\in L^{2n/(n+2)}(\R^n).
		\] 
		
		Considering one more time \(q\) near \(2^*_s\), it is clear that
		\[
			\int_{\R^n}(U_s)_t^q\,dx=t^{\frac{(n-2s)q}{2}-n}\int_{\R^n}U_s^q\,dx .
		\]
		Again, differentiating at \(t=1\),
		\[
			q\int_{\R^n}U_s^{q-1}\psi_s\,dx=\Big(\frac{(n-2s)q}{2}-n\Big)\int_{\R^n}U_s^q\,dx .
		\]
		This vanishes at \(q=2^*_s\). Differentiating in \(q\) at \(q=2^*_s\) we conclude
		\[
			2^*_s\int_{\R^n}\log(U_s)U_s^{p_s}\psi_s\,dx=\frac{n-2s}{2}\int_{\R^n}U_s^{2^*_s}\,dx.
		\]
	\end{proof}

	Let us now define \(K_s\coloneqq\{\phi\in\HH_{\rm rad}\colon\langle\phi,\psi_s\rangle_s=0\}\), and consider,
	for \(g\in\HH_{\rm rad}'\),
	\begin{equation}\label{eq:lin-s}
		\langle\phi,v\rangle_{\sigma,s}-\int_{\R^n}V_s\phi v\,dx=g(v)+c\int_{\R^n}V_s\psi_sv\,dx
		\quad\forall\,v\in\HH_{\rm rad},
	\end{equation}
	with unknowns \(\phi\in K_s\), \(c\in\R\).

	\begin{lemma}\label{orz-s}
		There exist \(\sigma_0,C>0\) such that for \(\sigma\in(0,\sigma_0)\) and
		\(g\in\HH_{\rm rad}'\), \eqref{eq:lin-s} has a unique solution, and
		\(\|\phi\|_{\sigma,s}+|c|\le C\|g\|_{\sigma,s,*}\).
	\end{lemma}
	\begin{proof}
		We follow the proof of Lemma~\ref{orz}. Assume \(\sigma_k\to0\),
		\(\|g_k\|_{\sigma_k,s,*}\to0\) and \(\|\phi_k\|_{\sigma_k,s}+|c_k|=1\). Testing with
		\(\psi_s\), since \(\langle\phi_k,\psi_s\rangle_s=\int_{\R^n} V_s\phi_k\psi_s\,dx=0\),
		\[
			c_k\int_{\R^n} V_s\psi_s^2\,dx=\sigma_k\int_{\R^n}\nabla\phi_k\cdot\nabla\psi_s\,dx-g_k(\psi_s),
		\]
		and
		\[
			\sigma_k\Big|\int_{\R^n}\nabla\phi_k\cdot\nabla\psi_s\,dx\Big|\le
			\sigma_k^{1/2}\|\phi_k\|_{\sigma_k,s}\|\nabla\psi_s\|_{L^2}\to0 ,
		\]
		so \(c_k\to0\). Up to a subsequence \(\phi_k\rightharpoonup\phi\) in \(\dot H^s(\R^n)\), and for
		radial \(v\in C^\infty_c\),
		\[
			\sigma_k\Big|\int_{\R^n}\nabla\phi_k\cdot\nabla v\,dx\Big|\le\sigma_k^{1/2}\|\nabla v\|_{L^2}\to0 .
		\]
		Hence \((-\Delta)^s\phi=V_s\phi\), \(\phi\) is radial and \(\langle\phi,\psi_s\rangle_s=0\), so
		\(\phi=0\) by nondegeneracy.

		We now justify the compactness used above: multiplication by
		\(V_s^{1/2}\in L^{n/s}(\R^n)\) is compact from \(\dot H^s(\R^n)\) into \(L^2(\R^n)\). We
		split the integral between a fixed ball \(B_R\) and its complement. On \(B_R\), the
		embedding \(\dot H^s(\R^n)\cap L^{2^*_s}(\R^n)\subset H^s_{\rm loc}(\R^n)\) together with
		Rellich's theorem gives \(\phi_k\to\phi\) strongly in \(L^2(B_R)\), and since
		\(V_s^{1/2}\) is bounded there, \(V_s^{1/2}\phi_k\to V_s^{1/2}\phi\) strongly in
		\(L^2(B_R)\) as well. Outside \(B_R\), H\"older's inequality bounds the contribution by
		\(\|V_s^{1/2}\|_{L^{n/s}(\{|x|>R\})}\), uniformly in \(k\) since \(\|\phi_k\|_{\sigma_k,s}\)
		is bounded, and this tail norm can be made as small as we like by taking \(R\) large,
		because \(V_s^{1/2}\in L^{n/s}(\R^n)\). Letting first \(k\to\infty\) and then
		\(R\to\infty\) gives \(V_s^{1/2}\phi_k\to V_s^{1/2}\phi=0\) in \(L^2(\R^n)\) (recall
		\(\phi=0\), as just shown), so
		\[
			\int_{\R^n} V_s\phi_k^2\,dx\to0
		\]
		and, testing with \(\phi_k\),
		\[
			\|\phi_k\|_{\sigma_k,s}^2=\int_{\R^n} V_s\phi_k^2\,dx+g_k(\phi_k)\to0 ,
		\]
		a contradiction. 
		
		Existence follows
		from the Fredholm alternative exactly as in Lemma~\ref{orz}.
	\end{proof}

	\subsection{Reduction}
	Fix \(\Lambda\coloneqq[\lambda_-,\lambda_+]\) with \(0<\lambda_-<\lambda_s<\lambda_+\) and let \(\sigma\coloneqq\lambda\eta\),
	\(\lambda\in\Lambda\). We look for \(z=U_s+\phi\), \(\phi\in K_s\), \(c\in\R\), with
	\begin{equation}\label{eq:redu-s}
		\langle\phi,v\rangle_{\sigma,s}-\int_{\R^n} V_s\phi v\,dx=\mathcal{E}_s(v)+\int_{\R^n} N_\eta(\phi)v\,dx+c\int_{\R^n} V_s\psi_sv\,dx
	\end{equation}
	for all \(v\in\HH_{\rm rad}\) where, using \((-\Delta)^sU_s=F_0(U_s)\),
	\[
		\mathcal{E}_s(v)\coloneqq\int_{\R^n}\left(F_\eta(U_s)-F_0(U_s)\right)v\,dx-\sigma\int_{\R^n}\nabla U_s\cdot\nabla v\,dx,
	\]
	\[
		N_\eta(\phi)\coloneqq\left[F_\eta(U_s+\phi)-F_\eta(U_s)-F_\eta'(U_s)\phi\right]
		+\left[F_\eta'(U_s)-F_0'(U_s)\right]\phi .
	\]

	\begin{lemma}\label{lemma:estimates-s}
		For \(\eta\) small, \(\lambda\in\Lambda\) and \(\phi,\phi_1,\phi_2\in\HH\),
		\begin{enumerate}[label=\textup{(\alph*)}]
			\item \(\|\mathcal{E}_s\|_{\sigma,s,*}\le C\eta\),
			\item \(\|N_\eta(\phi)\|_{\sigma,s,*}\le C\left(\|\phi\|_{\sigma,s}^{p_s+\eta}+\|\phi\|_{\sigma,s}^2+\eta\|\phi\|_{\sigma,s}\right)\),
			\item \(\begin{aligned}[t]
				\|N_\eta(\phi_1)&-N_\eta(\phi_2)\|_{\sigma,s,*}\le\\
				&\mkern-10mu C\left(\|\phi_1\|_{\sigma,s}^{p_s+\eta-1}+\|\phi_2\|_{\sigma,s}^{p_s+\eta-1}+\|\phi_1\|_{\sigma,s}+\|\phi_2\|_{\sigma,s}+\eta\right)\|\phi_1-\phi_2\|_{\sigma,s}.
			\end{aligned}\)
		\end{enumerate}
	\end{lemma}
	\begin{proof}
		As in Lemma~\ref{lemma:estimates}, we estimate \(\mathcal{E}_s\) and \(N_\eta\)
		separately. Since \(0<U_s\le1\), \(0\le1-U_s^\eta\le\eta|\log U_s|\) and
		\(|\log U_s|\le C\log(2+|x|)\). Thus
		\[
			|F_\eta(U_s)-F_0(U_s)|\le C\eta\log(2+|x|)(1+|x|)^{-(n+2s)},
		\]
		which has \(L^{q_s}\) norm \(O(\eta)\) because \((n+2s)q_s=2n\). Moreover
		\[
			\sigma\left|\int_{\R^n}\nabla U_s\cdot\nabla v\,dx\right|=\sigma\left|\int_{\R^n}(-\Delta U_s)v\,dx\right|
			\le\sigma\|\Delta U_s\|_{L^{q_s}}\|v\|_{L^{2^*_s}}
			\le C\eta \|v\|_{\sigma,s}
		\] 
		by \eqref{eq:Us-int} and \eqref{eq:sob-s}. For \(N_\eta,\) we use
		Lemma~\ref{lem:elem} with \(a=U_s\), \(p=p_s+\eta\), H\"older's inequality with
		\(\frac{p_s+\eta}{r_\eta}+\frac1{2^*_s}=1\), Lemma~\ref{lem:interp-s}, and, when \(p_s+\eta>2\)
		(which happens only if \(n\le6s\)), the bound
		\[
			\int_{\R^n} U_s^{p_s+\eta-2}\phi^2|v|\,dx\le\|U_s^{p_s+\eta-2}\|_{L^{w}}\|\phi\|_{L^{2^*_s}}^2\|v\|_{L^{2^*_s}}
		\]
		with
		\(\frac1w=1-\frac3{2^*_s}=\frac{6s-n}{2n}\) (and \(w=\infty\) if \(n=6s\)), here
		\(U_s^{p_s+\eta-2}\in L^w\) because \((n-2s)(p_s-2)w=2n>n\). We also use
		\(|F_\eta'(U_s)-F_0'(U_s)|\le C\eta(1+|\log U_s|)U_s^{p_s-1}\), whose \(L^{n/(2s)}\) norm is
		\(O(\eta)\).
	\end{proof}

	Exactly as in Lemma~\ref{eq:puntofijo}, Lemmas~\ref{orz-s} and
	\ref{lemma:estimates-s} and the contraction mapping principle give, for \(\eta\)
	small and \(\lambda\in\Lambda\), a unique solution \((\phi_{\eta,\lambda},c_{\eta,\lambda})\) of
	\eqref{eq:redu-s} with \(\|\phi_{\eta,\lambda}\|_{\sigma,s}\le C_0\eta\), depending continuously on
	\(\lambda\).

	\begin{lemma}\label{lemaimp-s}
		As \(\eta\to0\), uniformly in \(\lambda\in\Lambda\),
		\[
			c_{\eta,\lambda}\int_{\R^n} V_s\psi_s^2\,dx=\eta\left(-A_s+\lambda B_s+o(1)\right).
		\]
	\end{lemma}
	\begin{proof}
		Testing \eqref{eq:redu-s} with \(\psi_s\) and using
		\(\langle\phi,\psi_s\rangle_s=\int_{\R^n} V_s\phi\psi_s\,dx=0\),
		\begin{align*}
			c\int_{\R^n} V_s\psi_s^2\,dx&=\int_{\R^n}\left(F_0(U_s)-F_\eta(U_s)\right)\psi_s\,dx-\int_{\R^n} N_\eta(\phi)\psi_s\,dx\\
			&+\sigma B_s+\sigma\int_{\R^n}\nabla\phi\cdot\nabla\psi_s\,dx .
		\end{align*}
		As we did in Lemma \ref{lemaimp} we estimate each term on the right-hand side. Indeed we get that  \(F_0(U_s)-F_\eta(U_s)=c_{n,s}U_s^{p_s}(1-U_s^\eta)\) and
			\(|1-U_s^\eta+\eta\log U_s|\le\eta^2(\log U_s)^2\), so the first adding term is equals to
			\(-\eta A_s+O(\eta^2)\). The second term is bounded by
			\(\|N_\eta(\phi)\|_{\sigma,s,*}\|\psi_s\|_{\sigma,s}=O(\eta^{p_s}+\eta^2)=o(\eta)\).
			Finally \(\sigma B_s=\lambda\eta B_s\) and the last term is bounded by
			\(\sigma^{1/2}\|\phi\|_{\sigma,s}\|\nabla\psi_s\|_{L^2}=O(\eta^{3/2})\).
	\end{proof}

	\begin{proof}[Proof of Theorem~\ref{TeoremaBlowDown}]
		Since \(A_s,B_s>0\) and \(\lambda_s=A_s/B_s\) (by Lemma~\ref{taqui-s}, this is the value
		in the statement), \(-A_s+\lambda_-B_s<0<-A_s+\lambda_+B_s\). By Lemma~\ref{lemaimp-s} and
		continuity there is \(\lambda_\eta\in(\lambda_-,\lambda_+)\) with \(c_{\eta,\lambda_\eta}=0\), and
		\(\lambda_\eta\to\lambda_s\). Then \(z_\eta=U_s+\phi_{\eta,\lambda_\eta}\) is a weak solution of
		\eqref{elproblemaSigma} with \(\sigma_\eta=\lambda_\eta\eta\) (tested with radial
		functions, hence with all of \(\HH\) by Remark~\ref{rem:symcrit}). Its positivity
		follows from the proof of Lemma~\ref{lem:positivity}, applied to
		\[
			-\Delta z+\sigma_\eta^{-1}(-\Delta)^sz=\sigma_\eta^{-1}F_\eta(z),
		\]
		the only property of the
		right-hand side \(h=\sigma_\eta^{-1}F_\eta(z)\ge0\) used there is that \(h\) and \(h-h_k\) define
		nonnegative elements of \(\HH'\), which holds here because \(h\in L^{q_s}\) by
		Lemma~\ref{lem:interp-s}.

		Finally, \(\frac{2s}{p_s+\eta-1}=\frac{n-2s}{2}+O(\eta)\), so a change of variables gives, for
		every \(R>0\),
		\[
			\int_{B_R}u_\eta^{2^*_s}\,dy=\mu_\eta^{O(\eta)}\int_{B_{R/\mu_\eta}}z_\eta^{2^*_s}\,dx ,
		\]
		and \(\mu_\eta^{O(\eta)}\to1\) since \(\eta|\log\eta|\to0\). As \(z_\eta\to U_s\) in \(L^{2^*_s}\) and
		\(R/\mu_\eta\to0\), the right-hand side tends to \(0\) for \(R<\infty\) and to
		\(\int_{\R^n} U_s^{2^*_s}\,dx\) for \(R=\infty\).
	\end{proof}

\section*{Acknowledgments}
	The authors are grateful to the referees for their careful reading and
	valuable suggestions, which led to a substantial improvement of the paper.
	All the authors were partially supported by the project {\it An\'alisis de Fourier y Ecuaciones no locales en Derivadas Parciales}, Grant PID2023-148028NB-I00 funded by MICIU/AEI/10.13039/501100011033 and by ERDF/EU.
	AQ was also partially supported by FONDECYT Grant 1231585.
	LMDP was also supported by ANII  under the programs  “Movilizaciones AMSUD 2023-01”
	(MOV CO 2023 1 1012406 and MOV CO 9 101336) and ANII FCE 3 2024 1 181302.

	\bigskip

{\bf On behalf of all authors, the corresponding author states that there is no conflict of interest. }

	{\bf Data sharing is not applicable to this article, as no datasets were generated or analysed during the current study.}

\end{document}